\documentclass[12pt,a4paper]{article}
\usepackage[latin1]{inputenc}
\usepackage{amsmath}
\usepackage{amsthm}
\usepackage{amsfonts,dsfont}
\usepackage{amssymb}
\usepackage[dvipsnames]{xcolor}
\usepackage{hyperref,graphicx}

\newcommand{\vc}[1]{\boldsymbol{#1}}
\newcommand{\mbP}{\mathbb P}
\newcommand{\mbE}{\mathbb E}
\newcommand{\ee}{{\rm e}}
\newcommand{\dd}{{\rm d}}
\newcommand{\NN}{\mathbb{N}}
\newcommand{\eps}{\varepsilon}

\newtheorem{prop}{Proposition}
\newtheorem{theorem}{Theorem}
\newtheorem{lemma}{Lemma}
\newtheorem{cor}{Corollary}
\newtheorem{rem}{Remark}

\title{Consistent estimation in subcritical birth-and-death processes}
\author{Sophie Hautphenne\footnote{The University of Melbourne} \and Emma Horton\footnote{University of Warwick}}
\date{\today}

\begin{document}

\maketitle

\begin{abstract}

We investigate parameter estimation in subcritical continuous-time birth-and-death processes with multiple births. We show that the classical maximum likelihood estimators for the model parameters, based on the continuous observation of a single non-extinct trajectory, are not consistent in the usual sense: conditional on survival up to time $t$, they converge as $t \to \infty$ to the corresponding quantities in the associated $Q$-process, namely the process conditioned to survive in the distant future. We develop the first \emph{$C$-consistent} estimators in this setting, which converge to the true parameter values when conditioning on survival up to time $t$, and establish their asymptotic normality. The analysis relies on spine decompositions and coupling techniques.

\smallskip


\smallskip


\end{abstract}

\section{Introduction}\label{sec:intro}

Subcritical linear birth-and-death processes 
form a natural class of continuous-time Markov models for populations in decline. Such processes capture situations where each birth event can produce several offspring, but where, on average, the population fails to sustain itself and extinction eventually happens. This setting arises in many applied contexts. In ecology, subcritical dynamics can, for example, describe populations of endangered species that cannot persist without intervention 
\cite{Lande1993}. 
In epidemiology, subcritical branching processes model outbreaks that may cause clusters of secondary cases but cannot sustain long-term transmission~\cite{blumberg2013comparing}.
In cell biology, they can represent defective or pre-cancerous cell lineages that produce short bursts of offspring before dying out \cite{durrett2015branching,Nowak2006}. 
In all these scenarios, reliable statistical inference for the underlying birth and death rates, and for the offspring distribution, is essential for predicting extinction risk, assessing intervention strategies, and quantifying the expected persistence time of the population.

Despite this practical relevance, statistical inference for subcritical birth-and-death processes, based on the observation of a single long trajectory, has received little attention in the literature. The difficulty lies in the fact that these processes become extinct with probability one, which complicates inference under this observation scheme. Most of the literature instead considers the supercritical case, where populations have a positive probability of surviving indefinitely, and consistency of estimators is then studied under conditioning on survival of the trajectory \cite{Keiding1975,Guttorp1991}. See also \cite{hautphenne2024birth} for an overview of simulation and estimation methods for general birth-and-death processes. To the best of our knowledge, consistency of estimators from a single non-extinct trajectory of a continuous-time subcritical birth-and-death process has not previously been addressed.

In the subcritical case, the asymptotic $(t \to \infty)$ behaviour of the process conditioned to survive is quite different to that of the supercritical case. Indeed, the classical maximum likelihood estimators for the parameters of the birth-and-death process are not consistent in this setting: although they converge in probability as $t \to \infty$ when conditioning on $Z_t>0$, they do not converge to the true parameter values. This is due to the fact that, in the subcritical case, conditioning on survival alters the long-term behaviour of the process. As such, observations should be interpreted as being generated not by the original process, but by the so-called $Q$-process \cite{AthreyaNey1972,Lambert2007}, namely the process conditioned on ultimate survival.

In fact, in the subcritical case, a classical estimator $\hat{\theta}_t$ for a parameter $\theta$ of the process is only \emph{Q-consistent} \cite{BHM22}, in the sense that for any $\varepsilon > 0$,  
\begin{equation}\label{eq:Qconsistent}
\lim_{t \to \infty} \mathbb{P}\!\left[ \, \big| \hat{\theta}_t - \theta^\uparrow \big| > \varepsilon \,\middle|\, Z_t > 0 \right] = 0,
\end{equation}
where $\theta^\uparrow$ denotes the counterpart of $\theta$ in the associated $Q$-process.
 By contrast, we call an estimator $\tilde{\theta}_t$ \emph{C-consistent} if, for any $\varepsilon > 0$,  
\begin{equation}\label{eq:Cconsistent}
\lim_{t \to \infty} \mathbb{P}\!\left[ \, \big| \tilde{\theta}_t - \theta \big| > \varepsilon \,\middle|\, Z_t > 0 \right] = 0,
\end{equation}  
so that, conditional on survival up to time $t$, it converges to the true parameter value. This notion of consistency is particularly relevant in practice, since endangered populations are typically studied precisely because they are extant, and in that case accurate estimation of the model parameters is essential.

In this paper, we extend recent advances on conditional consistency in discrete-time settings, where $Q$-consistency and $C$-consistency were first formalised for Galton--Watson and population-size-dependent branching processes with almost sure extinction \cite{BHM22}, and $C$-consistent estimators were later developed for a class of parametric population-size-dependent branching processes \cite{BHM25}, to continuous-time subcritical birth-and-death processes with multiple births.  
Our main contribution is twofold: \emph{(i)} we show that the classical maximum likelihood estimators for the parameters of these processes are only $Q$-consistent, in the sense of \eqref{eq:Qconsistent}; and \emph{(ii)} we construct the first $C$-consistent estimators in this setting, in the sense of \eqref{eq:Cconsistent}, and establish their asymptotic normality.

The proofs rely on coupling techniques linking the process conditioned to survive until time~$t$ with the $Q$-process, which allow us to transfer asymptotic properties from the latter to the former. While related ideas were used in \cite{BHM22} to establish consistency results in the discrete-time setting, here we adapt and extend them to continuous-time birth-and-death processes with multiple births.

The paper is organised as follows. Section~\ref{sec:model} introduces the model and its spectral properties, and Section~\ref{sec:cond} describes the conditioned process and its relation to the $Q$-process. In Section~\ref{sec:C-con} we study consistent estimation of the birth and death rates $\lambda$ and $\mu$ under the assumption that the mean offspring $m$ is known: we show that the classical MLEs are $Q$-consistent, identify their asymptotic limits, and propose new estimators that are $C$-consistent and asymptotically normal. Section~\ref{sec:pk} extends these results to estimation of the offspring probabilities. Section~\ref{rem-est-m} discusses extensions to joint estimation of $\lambda, \mu,$ and $m$, and Section~\ref{sec:numerical} presents empirical analyses of the proposed estimators. Section~\ref{sec:proofs} collects the proofs and auxiliary lemmas, several of which rely on coupling techniques and spine decompositions. Appendix \ref{sec:appendix} provides details on the spectral properties of the process and the coupling between the process conditioned to survive until time $t$ and the $Q$-process.


\section{Subcritical birth-and-death processes}\label{sec:model}

Let $\NN := \{1, 2, \dots\}$ and let $\lambda, \mu > 0$ denote the birth and death rates, respectively. Let $(p_k)_{k \ge 2}$ be a probability distribution on $\{2, 3, \dots\}$, i.e., $p_k \in [0,1]$ for all $k \ge 2$ and $\sum_{k \ge 2} p_k = 1$.  

We consider a continuous-time birth-and-death process $Z = (Z_t)_{t \ge 0}$ taking values in $\NN \cup \{0\}$, defined as follows. Given $Z_0 = i \ge 1$, one of two events occurs:
\begin{itemize}
\item at rate $i\lambda$, the process jumps to $i-1+k$ with probability $p_k$;
\item at rate $i\mu$, the process jumps to $i-1$.
\end{itemize}
The first event corresponds to a birth: an individual is replaced by $k$ new individuals, or equivalently, one individual gives birth to $k-1$ offspring. The second event corresponds to a death, in which an individual is removed from the system. The state $0$ is absorbing: if $Z_t = 0$ for some $t \ge 0$, then $Z_{t+s} = 0$ for all $s \ge 0$.  
We write $\mbP_i$ for the law of $Z$ started from $i \ge 1$ individuals, and $\mbE_i$ for the corresponding expectation operator. We omit the subscript when $i=1$.  

Throughout the paper, let $\xi$ denote a random variable with distribution $(p_k)_{k \ge 2}$, called the \emph{offspring distribution}, i.e., $\mbP(\xi = k) = p_k$. We denote its mean and variance by
\[
m := \mbE[\xi], \qquad \sigma^2 := \mathrm{Var}(\xi),
\]
which we assume to be finite.

The mean growth rate of the process is
\begin{equation}\label{rho}
\rho := \lambda(m - 1) - \mu.
\end{equation}
To see this, define $\Psi_t := \mbE[Z_t]$ and condition on the first event (birth or death). Then
\[
\Psi_t = e^{-(\lambda + \mu)t} + \int_0^t \lambda e^{-(\lambda + \mu)u}\, m\, \Psi_{t-u}\, \mathrm{d}u, \qquad t \ge 0.
\]
Differentiating with respect to $t$ yields the ODE $\Psi_t' = \rho \Psi_t$ with initial condition $\Psi_0 = 1$, and hence $\Psi_t = e^{\rho t}$.  
Throughout the paper, we assume $\rho < 0$, i.e., the process is subcritical and therefore becomes extinct almost surely.

\section{The conditioned process}\label{sec:cond}

In this section, we study the process \( Z \) conditioned on survival up to time \( t \in (0, \infty] \), that is, conditioned on \( \{Z_t > 0\} \). Of particular interest is the asymptotic behaviour of \( Z \) conditional on $\{Z_t > 0\}$, as $t \to \infty$. As we will show, in this limit, the process admits a \emph{spine decomposition}, in which one distinguished individual (the spine) survives indefinitely and, at an accelerated rate, produces offspring according to a size-biased distribution, each initiating an independent copy of the original process.

Let $Q$ denote the sub-generator of $Z$ restricted to the transient states $\NN$, and let $P(t) = \exp(Q t)$ be the associated (sub-stochastic) transition semigroup for $t \geq 0$. One way to characterise the long-term behaviour of branching processes is via a Perron--Frobenius decomposition, which states that the leading-order behaviour of $P(t)$ is governed by the dominant eigentriple of $Q$. That is, there exists a constant $\rho_* \in \mathbb R$ and positive vectors $\vc{u}, \vc{v}$ such that
\begin{equation}\label{eigentriple}
  \vc u^\top Q = \rho_* \vc u^\top, \quad Q \vc v = \rho_* \vc v, \quad \vc u^\top \vc 1 = 1, \quad \vc u^\top \vc v = 1,
\end{equation}
and
\begin{equation}\label{PF-asymp}
  P(t) \sim e^{\rho_* t}\, \vc v \vc u^\top, \quad t \to \infty,
\end{equation}
where $\vc u^\top$ denotes the transpose of $\vc u$. The normalisation $\vc u^\top \vc 1 = 1$ implies that $\vc u$ is a probability distribution on $\mathbb N$ which, along with the normalisation $\vc u^\top \vc v = 1$, ensures that $\vc u$ and $\vc v$ are unique, i.e. not only unique up to multiplication by a constant. The existence of such a Perron--Frobenius triple in this setting is proved in Appendix~\ref{PFdec}.

An explicit calculation using \eqref{eigentriple} shows that $\rho_* = \rho$, where $\rho$ is given in \eqref{rho}, and that $v_j = c j$ for some normalising constant $c > 0$. The Perron--Frobenius eigenvectors have the following interpretation: $\vc u^\top$ corresponds to the \emph{quasi-stationary distribution} of $Z$, while $\vc v$ records the relative \emph{strength} or \emph{importance} of each state. Indeed, from \eqref{PF-asymp},
\[
\lim_{t\to\infty}\mbP_i(Z_t=j \mid Z_t>0)
  = \lim_{t\to\infty}\frac{\vc e_i^\top P(t) \vc e_j}{\vc e_i^\top P(t)\vc 1}
  = \lim_{t\to\infty}\frac{e^{\rho t} v_i u_j}{e^{\rho t} v_i} 
  = u_j, \qquad j \geq 1,
\]
and
\[
\lim_{t\to\infty}\frac{\mbP_j(Z_t>0)}{\mbP_i(Z_t>0)}
  = \lim_{t\to\infty}\frac{\vc e_j^\top P(t)\vc 1}{\vc e_i^\top P(t)\vc 1}
  = \frac{v_j}{v_i} = \frac{j}{i}, \qquad i,j \geq 1.
\]

We now consider the process $Z$ conditioned on survival up to time $t$, which defines a time-inhomogeneous Markov chain $(Z_u^{(t)})_{0 \le u \le t}$ with transition probabilities
\begin{align}
P^{(t)}_{ij}(u) 
&:= \mbP(Z_u^{(t)} = j \mid Z_0^{(t)} = i) \notag \\
&= \mbP(Z_u = j \mid Z_0 = i,\, Z_t > 0) \notag \\
&= P_{ij}(u)\, \frac{\vc e_j^\top P(t-u)\vc 1}{\vc e_i^\top P(t)\vc 1}, 
\qquad i,j \geq 1. \label{condProb}
\end{align}
Taking the limit as $t \to \infty$ in \eqref{condProb} and using the Perron--Frobenius asymptotic \eqref{PF-asymp}, we obtain the time-homogeneous transition probabilities
\begin{align}
P^\uparrow_{ij}(u) 
&:= \lim_{t \to \infty} P^{(t)}_{ij}(u) \notag \\
&= P_{ij}(u)\, \frac{v_j}{v_i}\, e^{-\rho u} \notag \\
&= P_{ij}(u)\, \frac{j}{i}\, e^{-[(m-1)\lambda - \mu]u}. \label{Pup}
\end{align}
Let $Z^\uparrow = (Z^\uparrow_t)_{t \geq 0}$ denote the Markov process with transition semigroup $P^\uparrow(t)$. When $\rho < 0$, this process is positive recurrent and is known in the literature as the \emph{$Q$-process} \cite{AthreyaNey1972}.

Let us now derive the generator $Q^\uparrow$ of $Z^\uparrow$ and give a probabilistic interpretation of the $Q$-process. From \eqref{Pup},
\begin{align*}
Q^\uparrow_{ij} 
&:= \lim_{t \to 0} \frac{d}{dt} P^\uparrow_{ij}(t) \\
&= \lim_{t \to 0} \frac{d}{dt}\Big( P_{ij}(t)\, \frac{j}{i}\, e^{-\rho t} \Big) \\
&= Q_{ij}\, \frac{j}{i} - \mathbf{1}_{\{j=i\}}\, \rho.
\end{align*}
Explicitly, we have
\begin{equation}\label{Qup}
Q^\uparrow_{ij} = 
\begin{cases}
(i-1)\mu, & j = i-1, \\[0.7em]
(i-1)\lambda p_k + \lambda m \cdot \frac{k p_k}{m}, & j = i-1+k, \; k \geq 2, \\[0.7em]
-(i-1)(\lambda+\mu) - m\lambda, & j = i.
\end{cases}
\end{equation}

This suggests the following interpretation of the $Q$-process. When initiated from $i$ individuals, one of them is chosen uniformly at random and marked the \emph{spine}. The remaining $i-1$ individuals each evolve as independent copies of the original process $(Z, \mbP)$. The spine gives birth at rate $m\lambda$, producing offspring according to the size-biased distribution $\tilde p_k := k p_k / m$. One of the $k$ offspring is selected uniformly at random to continue as the spine, while the other $k-1$ initiate independent copies of the original process, \cite{LPP95}. Alternatively, after removing the spine, the process $Z^\uparrow - 1$ can be viewed as a birth-and-death process with immigration, where births and deaths occur as in the original process $(Z,\mbP)$, while immigration occurs at rate $\lambda m$, with $\ell$ immigrants arriving with probability $\tilde{p}_{\ell+1}$ for $\ell \ge 1$. 

Another natural way to connect the measure $\mathbb P^\uparrow$ and the original measure $\mathbb P$ is via a martingale change of measure. Indeed, from the fact that $\Psi_t = {\rm e}^{\rho t}$, as shown at the end of the previous section, along with the Markov property, $({\rm e}^{-\rho t}Z_t)_{t \ge 0}$ defines a non-negative martingale, which can be used to define a new measure $\mathbb Q$ via 
\[
  \frac{{\rm d}\mathbb Q}{{\rm d}\mathbb P}\bigg|_{\sigma(Z_s, s \le t)} := {\rm e}^{-\rho t}Z_t. 
\]
By splitting on the first event, one can show that the dynamics of $Z$ under the measure $\mathbb Q$ agree with those of $(Z^\uparrow, \mathbb P^\uparrow)$.

This pathwise description yields the following useful identity: if $Z^\uparrow_0 = 1$ and $N_t$ denotes the number of births along the spine up to time $t$ (so that $N_t \sim \mathrm{Poi}(\lambda m t)$), with birth times $T_1,\dots,T_{N_t}$, then
\begin{equation}\label{spine_dec}
Z^\uparrow_t \;\overset{d}{=}\; 1 + \sum_{i=1}^{N_t} \;\sum_{\substack{j=1 \\ j \neq i^*}}^{\tilde{\xi}_i} Z^{(j)}_{t-T_i},
\end{equation}
where $i^*$ denotes the index of the spine, which is chosen uniformly from $\{1, \dots, \tilde\xi_i\}$, the $\{\tilde{\xi}_i\}$ are i.i.d.\ with size-biased offspring distribution $(\tilde p_k)_{k \geq 2}$, and the $Z^{(j)}$ are i.i.d.\ copies of $Z$ under $\mbP$, with $Z^{(j)}_{t-T_i}$ denoting the process initiated from the $j$-th individual born at time $T_i$.

We note that combining \eqref{PF-asymp} and \eqref{Pup} yields
\begin{equation}
  \lim_{t \to \infty} P^\uparrow_{ij}(t) = \lim_{t \to \infty}P_{ij}(t)\frac{v_j}{v_i}\ee^{-\rho t} = u_jv_j, \qquad i, j \ge 1.
\label{Q-st-dist}
\end{equation}
In other words, because $v_j=cj$, the stationary distribution of the $Q$-process is the vector whose $j$-th entry is given by $j\,u_j/(\sum_{j} j\,u_j)$.
Let $$\pi^\uparrow:=\sum_j j \, u_j v_j= \frac{\sum_j j^2 \, u_j}{\sum_{j} j\,u_j} $$ denote the mean of this stationary distribution. This quantity will play an important role throughout the paper. In the next lemma we derive a closed-form expression for $\pi^\uparrow$. Since the proof relies on a direct application of the spine decomposition \eqref{spine_dec}, we present it here rather than postponing it to Section~\ref{sec:proofs}.
\begin{lemma}\label{lem:pi-up}
We have
\[
\pi^\uparrow = 1 - \frac{\lambda[\sigma^2 + m(m-1)]}{\rho}.
\]
\end{lemma} 

\noindent\textbf{Proof.}
From \eqref{Q-st-dist}, 
\[
\pi^\uparrow = \sum_{j \ge 1} j \, u_j v_j = \lim_{t \to \infty} \mathbb E_i[Z_t^\uparrow], \quad \text{for any } i \ge 1.
\]
It suffices to compute this limit for $i=1$, since for $i \ge 1$ one has
\begin{equation}\label{ab}
  \mathbb E_i[Z_t^\uparrow] = \mathbb E_1[Z_t^\uparrow] + \mathbb E_{i - 1}[Z_t],
\end{equation}
and, because $Z$ is subcritical ($\rho<0$), $\mathbb E_{i-1}[Z_t] = (i-1) e^{\rho t} \to 0$ as $t \to \infty$.

From the spine decomposition \eqref{spine_dec}, conditioning on the values $k \ge 2$ of $\tilde{\xi}$, and using the fact that one of the $k$ offspring is chosen uniformly to continue the spine while the remaining $k-1$ initiate independent copies of the original process $(Z,\mbP)$, we obtain
\[
\mbE_1[Z_t^\uparrow] 
= 1 + \mbE\bigg[\sum_{i = 1}^{N_t}\sum_{k \ge 2}\frac{k p_k}{m} \sum_{i^* = 1}^k \frac1k \sum_{\substack{j = 1 \\ j \neq i^*}}^k Z_{t-T_i}^{(j)} \bigg].
\]
Since $N_t \sim \mathrm{Poi}(\lambda m t)$ and, conditional on $N_t$, the birth times are i.i.d.\ uniform on $[0,t]$, it follows that
\begin{align}
\mbE_1[Z_t^\uparrow] 
&= 1 + \frac{\lambda m t}{t}\sum_{k \ge 2} \frac{k(k-1)p_k}{m}\int_0^t \mbE[Z_{t-s}]\,\dd s \notag\\
&= 1 + \lambda\big(\sigma^2 + m(m-1)\big)\int_0^t e^{\rho s}\,\dd s \notag\\
&= 1 + \lambda\big(\sigma^2 + m(m-1)\big)\frac{e^{\rho t}-1}{\rho}. \label{M1-up}
\end{align}
Letting $t \to \infty$ gives
\[
\pi^\uparrow = \lim_{t \to \infty} \mbE_1[Z_t^\uparrow] 
= 1 - \frac{\lambda\big(\sigma^2 + m(m-1)\big)}{\rho},
\]
as claimed. \hfill $\square$

\section{Consistent estimation of the birth and death rates}\label{sec:C-con}

We are interested in estimating the birth and death rates, $\lambda$ and $\mu$, from the continuous observation of a population over the time interval $[0,t]$ for some $t>0$. We let
\begin{itemize}
  \item $b_t$ denote the total number of birth \textit{events} in $[0,t]$ (note that this differs from the total number of offspring),
  \item $d_t$ denote the total number of death events in $[0,t]$, and
  \item $\tau_t := \int_0^t Z_s \,{\rm d} s$ denote the total cumulative lifetime of all individuals over $[0,t]$.
\end{itemize}
We first assume that the offspring mean, $m$, is known; we discuss the case where it is unknown in Section \ref{rem-est-m}.

It is well known that in the supercritical case (i.e.\ when \( \rho > 0 \)), the classical maximum likelihood estimators (MLEs)
\begin{equation}\label{MLEs}
\hat{\lambda}_t := \frac{b_t}{\tau_t}, \qquad \hat{\mu}_t := \frac{d_t}{\tau_t},
\end{equation}
are consistent for \( \lambda \) and \( \mu \), respectively. That is, \( \hat{\lambda}_t \to \lambda \) and \( \hat{\mu}_t \to \mu \) in probability as \( t \to \infty \) on the set of non-extinction; see \cite{athreya1977estimation}.
This means that they are \( C \)-consistent estimators, where we recall that an estimator \( \hat\theta_t \) is called \( C \)-consistent for a quantity \( \theta \) associated with \( Z \) if, for any \( i \ge 1 \) and any \( \varepsilon > 0 \),
\begin{equation}\label{eq:Cconsistent2}
  \lim_{t \to \infty} \mathbb{P}_i(|\hat\theta_t - \theta| > \varepsilon \mid Z_t > 0) = 0. 
\end{equation}

However, \( C \)-consistency of the MLEs does \emph{not} hold in the subcritical setting. The heuristic reason is that when observing subcritical populations, we only see those individuals that have not yet become extinct. Thus, we should interpret the data as arising from the process conditioned on \( \{Z_t > 0\} \), which naturally introduces a bias into the estimates. For long (non-extinct) trajectories, the observed data can be interpreted as coming from the \( Q \)-process, rather than from the original process.

In the next proposition, we characterise the conditional limits of the MLEs \( \hat{\lambda}_t \) and \( \hat{\mu}_t \) given \( \{Z_t > 0\} \), in the subcritical regime.

\begin{prop}[$Q$-consistency of \( \hat{\lambda}_t \) and \( \hat{\mu}_t \)] \label{prop:Qcon}
In the subcritical case \( \rho < 0 \), for any \( i \ge 1 \) and \( \varepsilon > 0 \),
\begin{equation}
\lim_{t \to \infty} \mathbb{P}_i\left(|\hat{\lambda}_t - \lambda^\uparrow| > \varepsilon \mid Z_t > 0\right) = 0,
\label{Q-con-lam}
\end{equation}
and
\begin{equation}
\lim_{t \to \infty} \mathbb{P}_i\left(|\hat{\mu}_t - \mu^\uparrow| > \varepsilon \mid Z_t > 0\right) = 0,
\label{Q-con-mu}
\end{equation}
where
\begin{equation}\label{lam_mu_uparrow}
  \lambda^\uparrow := \frac{\lambda(\pi^\uparrow - 1) + \lambda m}{\pi^\uparrow}, \quad\text{ and }\quad 
  \mu^\uparrow := \frac{\mu (\pi^\uparrow - 1)}{\pi^\uparrow}.
\end{equation}
\end{prop}

\begin{rem}\label{rem_int}
\begin{enumerate}
\item
The limits $\lambda^\uparrow$ and $\mu^\uparrow$ can be interpreted as the analogous quantities to $\lambda$ and $\mu$ in the $Q$-process. To see this, recall that the $Q$-process consists of a single immortal individual (the spine) giving birth at rate $\lambda m$ to independent copies of the original process. Since $\pi^\uparrow$ represents the long-run average population size in the $Q$-process, we may interpret this as $\pi^\uparrow - 1$ individuals evolving as in the original process, each giving birth at rate $\lambda$, together with one spine individual giving birth at rate $\lambda m$. The numerator of $\lambda^\uparrow$ thus corresponds to the total long-run birth rate across the population, and dividing by $\pi^\uparrow$ yields the effective per-individual birth rate in the $Q$-process. A similar interpretation holds for $\mu^\uparrow$.

\item Observe that, since $m\geq 2$, we have \( \lambda^\uparrow\geq \lambda \), and \( \mu^\uparrow\leq \mu \), which is consistent with the fact that the $Q$-process does not become extinct, whereas the original process does.

\item
In the deterministic case where $p_m = 1$ for some $m \ge 2$ (that is, each birth event produces exactly $m$ children and $\sigma^2=0$), 
we obtain
\[
\pi^\uparrow = \frac{\mu + \lambda (m-1)^2}{\mu - \lambda (m-1)},
\]
and the following direct relationship holds between $\lambda^\uparrow$ and $\mu^\uparrow$:
\[
\lambda^\uparrow(m-1) - \mu^\uparrow = 0.
\]
In particular, in the binary case ($m=2$),
\[
\lambda^\uparrow = \mu^\uparrow = \frac{2\lambda\mu}{\lambda+\mu}.
\]

\end{enumerate}
\end{rem}

We can formalise the intuition given in Remark \ref{rem_int} by showing that $\lambda^\uparrow$ and $\mu^\uparrow$ are the respective limits of the analogues of the classical MLEs in the $Q$-process. That is, we consider
\begin{equation}\label{MLEsQ}
\hat{\lambda}^\uparrow_t := \frac{b_t^\uparrow}{\tau_t^\uparrow}, \qquad \hat{\mu}_t := \frac{d_t^\uparrow}{\tau_t^\uparrow},
\end{equation}
where $b_t^\uparrow$, $d_t^\uparrow$, and $\tau_t^\uparrow$ denote, respectively, the total number of birth events, death events, and the total cumulative lifetimes of individuals in the time interval $[0,t]$ in the $Q$-process. We then have the following result.

\begin{prop}[Consistency of $\hat{\lambda}^\uparrow_t$ and $\hat{\mu}^\uparrow_t$]\label{prop:Con_Qprocess}
For any $i \ge 1$ and $\varepsilon>0$,
\begin{equation}
\lim_{t\to\infty} \mbP_i\left(|\hat{\lambda}^\uparrow_t - \lambda^\uparrow| > \varepsilon\right) = 0,
\label{Q-con-lam2}
\end{equation}
and
\begin{equation}
\lim_{t\to\infty} \mbP_i\left(|\hat{\mu}^\uparrow_t - \mu^\uparrow| > \varepsilon\right) = 0.
\label{Q-con-mu2}
\end{equation}
\end{prop}

Propositions \ref{prop:Qcon} and \ref{prop:Con_Qprocess} indicate that, in the subcritical case, the MLEs $\hat{\lambda}_t$ and $\hat{\mu}_t$ satisfy a  different notion of consistency, known as $Q$-consistency. We say that an estimator $\hat\theta_t$ is $Q$-consistent for a quantity $\theta$ associated with the process $Z$ if, for all $i \ge 1$ and $\eps > 0$,
\begin{equation}\label{eq:Qconsistent2}
  \lim_{t \to \infty}\mathbb P_i\left(|\hat\theta_t - \theta^\uparrow| > \eps \mid Z_t > 0\right) = 0,
\end{equation}
where $\theta^\uparrow$ is the analogue of $\theta$ for the $Q$-process, $Z^\uparrow$; see \cite[Definition 1]{BHM22}.

As previously stated, our aim is to construct $C$-consistent estimators for $\lambda$ and $\mu$. To do this, we must correct for the bias introduced by the dynamics of the spine in the observed data.
To this end, we define the adjusted estimators
\begin{equation}\label{C-estimators}
\tilde\lambda_t := \frac{b_t}{\tau_t + (m-1)t}=\frac{\tau_t }{\tau_t +(m-1)t}\,\hat{\lambda}_t, \qquad \tilde\mu_t := \frac{d_t}{\tau_t - t}=\frac{\tau_t }{\tau_t -t}\,\hat{\mu}_t.
\end{equation}
To interpret these estimators, recall that in the interval $[0, t]$, the spine lives for $t$ units of time and produces offspring at rate $m\lambda$, while the remaining individuals reproduce at rate $\lambda$. Thus, the expected total number of birth events until time $t$ in the $Q$-process is 
$$\mbE [b_t^\uparrow]=\lambda (\mbE[\tau_t^\uparrow]-t)+\lambda m t=\lambda \left(\mbE[\tau_t^\uparrow]+(m-1)t\right),$$ which provides an intuitive interpretation for $\tilde\lambda_t$. Similarly, the expected total number of death events until time $t$ in the $Q$-process is 
$$\mbE [d_t^\uparrow]=\mu (\mbE[\tau_t^\uparrow]-t),$$ that is,
since the spine is immortal, we subtract its contribution $t$ from the cumulative lifetimes $\tau_t$ to estimate  $\mu$. This leads to the following result.

\begin{theorem}[$C$-consistency of $\tilde{\lambda}_t$ and $\tilde{\mu}_t$]\label{thm:C-consistency}
The estimators $\tilde\lambda_t$ and $\tilde\mu_t$ are $C$-consistent for $\lambda$ and $\mu$, respectively. 
That is, for any $i \ge 1$ and $\eps > 0$, 
\begin{equation}
\lim_{t\to\infty} \mbP_i(|\tilde{\lambda}_t - \lambda |>\varepsilon\,|\,Z_t>0)=0,
\label{C-con-lam}
\end{equation}
and
\begin{equation}
\lim_{t\to\infty} \mbP_i(|\tilde{\mu}_t -\mu |>\varepsilon\,|\,Z_t>0)=0.
\label{C-con-mu}
\end{equation}
\end{theorem}

We now establish the asymptotic distribution of $\tilde{\lambda}_t$ and $\tilde{\mu}_t$.

\begin{theorem}[Asymptotic normality of $\tilde{\lambda}_t$ and $\tilde{\mu}_t$]\label{thm:CLT}
Conditional on $\{Z_t > 0\}$, we have
\begin{equation}
  \frac{\sqrt{t}\,(\tilde\lambda_t - \lambda)}{\sqrt{\lambda/(\pi^\uparrow + m - 1)}} \xrightarrow{d} Y_1,
\label{CLT-lam}
\end{equation}
and
\begin{equation}
  \frac{\sqrt{t}\,(\tilde\mu_t - \mu)}{\sqrt{\mu/(\pi^\uparrow - 1)}} \xrightarrow{d} Y_2,
\label{CLT-mu}
\end{equation}
where $Y_1$ and $Y_2$ are independent standard normal random variables. 
\end{theorem}

The next corollary highlights that the estimator for $\lambda$ is more efficient than that for $\mu$. The difference arises from the fact that despite conditioning on survival, the original process is subcritical and so the overall effect results in a process that survives but the population size remains fairly small. As such, most of the time, birth events occur just often enough to keep the population from going extinct, with very rare deviations from this behaviour. However, the number of death events can vary much more while still remaining in this regime.

\begin{cor}\label{cor-asym-var}
Conditional on $\{Z_t > 0\}$, as $t \to \infty$, we have
\[
  {\rm Var}(\tilde\lambda_t) < {\rm Var}(\tilde\mu_t). 
\]
\end{cor}

\section{Consistent estimation of the offspring distribution}\label{sec:pk}

In this section, we discuss the consistent estimation of the probability $p_k$ that, at a birth event, the parent individual splits into $k$ new individuals, for a given $k\geq 2$. 

We note that, unless one assumes a parametric form for the offspring distribution or restricts the support of $(p_k)_{k\geq 2}$ to be finite or truncated, it is not possible to estimate the entire distribution $(p_k)_{k\geq 2}$ at once. Here we focus on the estimation of a single component $p_k$ at a time, treating the other probabilities as fixed or unknown, and without enforcing the normalisation constraint $\sum_k p_k=1$. This leads to a simple closed-form $C$-consistent estimator for $p_k$, which has a direct interpretation in terms of the $Q$-process and provides a tractable framework for asymptotic analysis. We discuss the joint estimation of $(p_2,p_3,\dots)$ in the finite support case in Remark~\ref{rem-multinomial} and Section~\ref{rem-est-m}.

For each $k\geq 2$, the classic MLE for $p_k$ on the set $\{b_t>0\}$ is
\begin{equation}
\hat{p}_{k,t}:= \dfrac{b_{k,t}}{b_t},
\end{equation}
where  $b_{k,t}$ denotes the total number of birth {events} that generate $k$ offspring in the time interval $[0,t]$, and $b_t$ is the total number of birth events. This estimator is $C$-consistent in the supercritical case, but not in the subcritical case (see Corollary \ref{prop:Qcon_pk}). 

In the subcritical case, we propose the following $C$-consistent estimator:
\begin{equation}\label{hoho}
\tilde{p}_{k,t}:= \dfrac{b_{k,t}}{{\lambda}\,[\tau_t+(k-1)\,t]}=\dfrac{b_{t}}{{\lambda}\,[\tau_t+(k-1)\,t]}\,\hat{p}_{k,t}.
\end{equation}
This is  motivated by the fact that the expected number of birth events generating $k$ offspring by time $t$ in the $Q$-process satisfies
$$\mbE [b_{k,t}^\uparrow]=\lambda \,p_k\, (\mbE[\tau_t^\uparrow]-t)+\lambda m \left(\frac{k\, p_k}{m}\right)\,t=\lambda\, p_k \left[\mbE[\tau_t^\uparrow]+(k-1)t\right].$$
In practice, if $\lambda$ is unknown, we replace it by the $C$-consistent estimator $\tilde{\lambda}_t$, so that $\tilde p_{k, t}$ would then be equal to
\begin{equation*}
\dfrac{b_{k,t}}{\tilde{\lambda}_t\,(\tau_t+(k-1)\,t)}=\dfrac{\tau_t+(m-1)\,t}{\tau_t+(k-1)\,t}\,\hat{p}_{k,t}.
\end{equation*}

\begin{theorem}[$C$-consistency and asymptotic normality of $\tilde{p}_{k,t}$]\label{thm:C-consistency_pk}
For any $k\geq 2$, the estimator $\tilde{p}_{k,t}$ is $C$-consistent for $p_k$, 
that is, for any $i \ge 1$ and $\eps > 0$, 
\begin{equation}
\lim_{t\to\infty} \mbP_i(|\tilde{p}_{k,t} - p_k |>\varepsilon\,|\,Z_t>0)=0.
\label{C-con-pk}
\end{equation}
Furthermore, conditional on $Z_t>0$, we have
\begin{equation}
  \frac{\sqrt{t}\,(\tilde{p}_{k,t} - p_k)}{\sqrt{p_k/(\lambda(\pi^\uparrow + k - 1))}} \xrightarrow{d} Y,
\label{CLT-pk}
\end{equation}
where $Y$ is a standard normal random variable.
\end{theorem}

We are  able to identify the conditional limit of the classical MLE $\hat{p}_{k,t}$ in the subcritical case.
\begin{cor}
[$Q$-consistency of \( \hat{p}_{k,t} \)] \label{prop:Qcon_pk}
In the subcritical case \( \rho < 0 \), for any \( i \ge 1 \) and \( \varepsilon > 0 \),
\begin{equation}
\lim_{t \to \infty} \mathbb{P}_i\left(|\hat{p}_{k,t} - p_k^\uparrow| > \varepsilon \mid Z_t > 0\right) = 0,
\label{Q-con-pk}
\end{equation}
where
\begin{equation}\label{pk_uparrow}
  p_k^\uparrow := p_k \,\frac{\pi^\uparrow +k- 1 }{\pi^\uparrow+m-1}.
\end{equation}
\end{cor}

\begin{rem}
Similar to the conditional limits $\lambda^\uparrow$ and $\mu^\uparrow$ of the classical MLEs for $\lambda$ and $\mu$ in a subcritical birth-and-death process (see Remark \ref{rem_int}), the limit $p_k^\uparrow$ can be interpreted as the counterpart of $p_k$ in the $Q$-process. Indeed, $p_k^\uparrow$ can be written as
$$p_k^\uparrow=\dfrac{\lambda\, p_k\,(\pi^\uparrow-1)+(\lambda \,m)(k\,p_k/m) }{\sum_{k\geq 2}\left[\lambda\, p_k\,(\pi^\uparrow-1)+(\lambda \,m)(k\,p_k/m)\right] },$$
where the numerator represents the asymptotic rate at which a birth event resulting in $k$ offspring occurs, and the denominator represents the asymptotic total birth rate, in the $Q$-process.

\end{rem}

\begin{rem}\label{rem-multinomial}
The estimator \eqref{hoho} can be viewed as the analogue of the `Poisson MLE', obtained by treating each $p_k$ in isolation without enforcing the normalisation constraint. If the offspring distribution is assumed to have finite support $\{2,\dots,M\}$, we may  enforce $\sum_{k=2}^M p_k=1$ and jointly derive 
\[
{\bar{p}}_{k,t}
=\dfrac{b_{k,t}/[\tau_t+(k-1)\,t]}{\sum_{\ell=2}^M b_{\ell,t}/[\tau_t+(\ell-1)\,t]},\qquad k=2,\dots,M,
\]
which is the analogue of the `multinomial MLE'. Both approaches yield $C$-consistent estimators, and the two coincide asymptotically since the normalising constant in the denominator of ${\bar{p}}_{k,t}$ tends to~1. We do not derive the asymptotic properties of the multinomial version here, as the presence of the normalisation constraint makes the analysis substantially more involved. In the classical ($Q$-consistent) setting, the same distinction arises, and the constrained estimator converges to the same limit, $p_k^\uparrow $, as the unconstrained estimator, ${\hat{p}}_{k,t}$.
\end{rem}

\section{Consistent joint estimation of $\lambda$, $\mu$ and $m$}\label{rem-est-m}

In this section, we discuss several approaches to the joint estimation of $\lambda$, $\mu$, and $m$ in a $C$-consistent framework. Our goal is to outline the main ideas behind these approaches, without entering into a formal analysis of their asymptotic properties. Instead, their asymptotic behaviour is examined empirically in Section~\ref{sec:numerical}.

\subsection{The $\delta$-skeleton approach}\label{delt}

One way to jointly estimate several parameters is to work with discrete-time statistics obtained from observing skeletons of the process. For $\delta > 0$, the $\delta$-skeleton of $Z$ is the discrete-time process $(Z_{n\delta})_{n \ge 0}$ obtained by observing $Z$ at multiples of $\delta$.
This is an embedded Galton-Watson (GW) process whose offspring mean and variance are given by 
\begin{align}
m^* &= e^{[\lambda(m - 1) - \mu] \delta}, \label{me1} \\
\sigma^{2*} &= e^{[\lambda(m - 1) - \mu] \delta} \left[1 - e^{[\lambda(m - 1) - \mu] \delta} \right] \,\pi^\uparrow,
 \label{me2}
\end{align}
where recall from Lemma \ref{lem:pi-up} that $\pi^\uparrow$ is a function of the mean $m$ and variance $\sigma^2$ of the original offspring distribution.
These expressions can be derived from the backward Kolmogorov differential equation satisfied by the p.g.f. \( F(t,s) \) of the population size \( Z_t \) starting from one individual at time \( t = 0 \):
$$\dfrac{\partial F(t,s)}{\partial t}=\mu-(\lambda+\mu)\,F(t,s)+\lambda\,P(F(t,s)),\quad F(0,s)=s,$$ where $P(s):=\sum_{k\geq 2} p_k s^k$ is the p.g.f.\ of the offspring distribution.

Proposition 4 in \cite{BHM22} provides \( C \)-consistent estimators \( \hat{m}^*_n \) and \( \hat{\sigma}^{2*}_n \) for $m^*$ and $\sigma^{2*}$, based on a least squares approach. The two moment equations \eqref{me1}--\eqref{me2}, together with the two additional  equations from \eqref{lam_mu_uparrow}, form a system of four equations in the four unknowns \( (\lambda, m, \mu, \sigma^2) \). Assuming this system is locally identifiable, we can, by the Continuous Mapping Theorem, jointly construct \( C \)-consistent estimators for all four parameters.
 
 In the binary case where \( p_2 = 1 \), we have \( m = 2 \) and \( \sigma^2 = 0 \), and the expressions for the mean and variance in the \(\delta\)-skeleton simplify to
\begin{align}
m^* &= \exp\left\{ (\lambda - \mu)\delta \right\}, \\
\sigma^{2*} &= \frac{\lambda + \mu}{\lambda - \mu} \exp\left\{(\lambda - \mu)\delta\right\} \left( \exp\left\{(\lambda - \mu)\delta\right\} - 1 \right).
\end{align}
The \( C \)-consistent estimators \( \hat{m}^*_n \) and \( \hat{\sigma}^{2*}_n \) from \cite[Proposition 4]{BHM22} then yield the following \( C \)-consistent estimators for \( \lambda \) and \( \mu \):
\begin{align}
\hat{\lambda}^*_n &= \frac{\log(\hat{m}^*_n)}{2\delta} \left( \frac{\hat{\sigma}^{2*}_n}{\hat{m}^*_n(\hat{m}^*_n - 1)} + 1 \right), \\
\hat{\mu}^*_n &= \frac{\log(\hat{m}^*_n)}{2\delta} \left( \frac{\hat{\sigma}^{2*}_n}{\hat{m}^*_n(\hat{m}^*_n - 1)} - 1 \right).
\end{align}

 While the $\delta$-skeleton approach may be effective, particularly as \( \delta \to 0 \), it may suffer from greater statistical variance than the estimators proposed in the previous sections (see Section~\ref{sec:numerical}). This is to be expected, as our new estimators are of the MLE type, which are known to be asymptotically efficient.

\subsection{The $Q$-process MLE approach}\label{Qpro}

Another approach consists of assuming that our observations come directly from the $Q$-process $Z^{\uparrow}$ and jointly estimating $\lambda$, $\mu$, and $m$ as parameters of $Z^{\uparrow}$ via maximum likelihood. This is motivated by the fact that the counterpart of the $C$-consistent estimators for $\lambda$, $\mu$, and $p_k$ in $Z^{\uparrow}$ are the MLEs for the corresponding quantities in $Z^{\uparrow}$ (see \eqref{tilde_up} and \eqref{hihi}).

From \eqref{Qup}, we note that the $Q$-process is a birth-and-death process with total birth and death rates at population size $r$ given by
$$\lambda_r = \sum_{\ell\geq 2}\lambda\,p_\ell\,(r+\ell-1)=\lambda\,(r-1+m), \quad\text{and}\quad \mu_r = \mu\,(r-1).$$
Assuming that a trajectory of $Z^{\uparrow}$ is continuously observed over the interval $[0,t]$, and the observations are recorded in $\vc X$, the log-likelihood function is given by
\begin{align}
\ell(\vc X, \lambda, \mu, m)
&= -\sum_{r \geq 1} [\lambda\,(r - 1 + m) + \mu\,(r - 1)]\,\nu_{r,t}^\uparrow \nonumber \\
&\quad + \sum_{r \geq 1} \beta_{r,t}^\uparrow \log[\lambda\,(r - 1 + m)] + \sum_{r \geq 2} \delta_{r,t}^\uparrow \log[\mu\,(r - 1)],
\label{LL2}
\end{align}
where $\nu_{r,t}^\uparrow$ denotes the total time spent in state $r$ during $[0,t]$, $\beta_{r,t}^\uparrow$ the number of birth events in state $r$, and $\delta_{r,t}^\uparrow$ the number of death events.

The MLEs $\tilde{\lambda}_t^\uparrow$, $\tilde{\mu}_t^\uparrow$, and $\tilde{m}_t^\uparrow$ of $\lambda$, $\mu$ and $m$ are solutions of the corresponding score equations. Using the identities $\sum_{r \geq 1} r\,\nu_{r,t}^\uparrow = \tau_t^\uparrow$, $\sum_{r \geq 1} \nu_{r,t}^\uparrow = t$, $\sum_{r \geq 1} \beta_{r,t}^\uparrow = b_t^\uparrow$, and $\sum_{r \geq 1} \delta_{r,t}^\uparrow = d_t^\uparrow$, we obtain the following system to solve in 
$\tilde{\lambda}_t^\uparrow$, $\tilde{\mu}_t^\uparrow$, and $\tilde{m}_t^\uparrow$:
\begin{align*}
\tilde{\lambda}_t^\uparrow &= \dfrac{b_t^\uparrow}{\tau_t^\uparrow + (\tilde{m}_t^\uparrow - 1)\,t}, \\
\tilde{\mu}_t^\uparrow &= \dfrac{d_t^\uparrow}{\tau_t^\uparrow - t}, \\
\tilde{\lambda}_t^\uparrow \,t &= \sum_{r \geq 1} \dfrac{b_r^\uparrow}{r - 1 + \tilde{m}_t^\uparrow}.
\end{align*}
We observe that $\tilde{\mu}_t^\uparrow$ has a closed-form expression (as in \eqref{tilde_up}), and is decoupled from the other equations. Substituting the expression for $\tilde{\lambda}_t^\uparrow$ into the last equation yields an implicit equation for $\tilde{m}_t^\uparrow$:
\[
\dfrac{b_t^\uparrow}{\tau_t^\uparrow/t + (\tilde{m}_t^\uparrow - 1)} = \sum_{r \geq 1} \dfrac{b_r^\uparrow}{r - 1 + \tilde{m}_t^\uparrow},
\]
which does not admit a closed-form solution. In practice, if $R_t$ denotes the maximum population size observed over $[0,t]$, $\tilde{m}_t^\uparrow$ satisfies the fixed-point equation $x=\mathcal{F}(x)$ where $$\mathcal{F}(x):=1- \dfrac{\tau_t^\uparrow}{t}+\dfrac{b_t^\uparrow}{\sum_{r = 1}^{R_t} \frac{b_r^\uparrow}{(r - 1 + x)}},$$which can be solved numerically.

In the case where the offspring distribution has finite support $\{2,3,\ldots,M\}$, a related approach consists in jointly estimating $\lambda$, $\mu$, and $p_2, \ldots, p_{M}$ by imposing the constraint $\sum_{k=2}^M p_k=1$, and using the fact that $m = \sum_{\ell=2}^M \ell\,p_\ell$. Skipping algebraic details, the resulting $C$-consistent estimators are explicitly given by
 \[
{\bar{p}}_{k,t}^\uparrow=\frac{b_{k,t}^\uparrow/[\tau_t^\uparrow+(k-1)t]}{\sum_{\ell=2}^M b_{\ell,t}^\uparrow/[\tau_t^\uparrow+(\ell-1)t]},\qquad 2\le k\le M,\]
\[
{\bar{\lambda}}_t^\uparrow=\sum_{k=2}^M \frac{b_{k,t}^\uparrow}{\tau_t^\uparrow+(k-1)t},\qquad
\tilde\mu_t^\uparrow=\frac{d_t^\uparrow}{\tau_t^\uparrow-t},
\]
where here $b_{k,t}^\uparrow$ denotes the number of births of size $k$ during [0, t] (see also Remark \ref{rem-multinomial} on joint estimators for $p_k$).

We note that the above assumes a non-parametric offspring distribution. If the offspring distribution belongs to a parametric family, the parameters can similarly be estimated jointly with $\lambda$ and $\mu$ via maximum likelihood.

\section{Numerical illustrations}\label{sec:numerical}

Through simple examples, we demonstrate the substantial bias of the classical MLEs for the parameters of a subcritical birth-and-death process when based on a single trajectory. This motivates the use of $C$-consistent estimators, which correct for the bias.

\subsubsection*{Binary case}

\begin{figure}[h!]
\centering\includegraphics[width=0.6\textwidth]{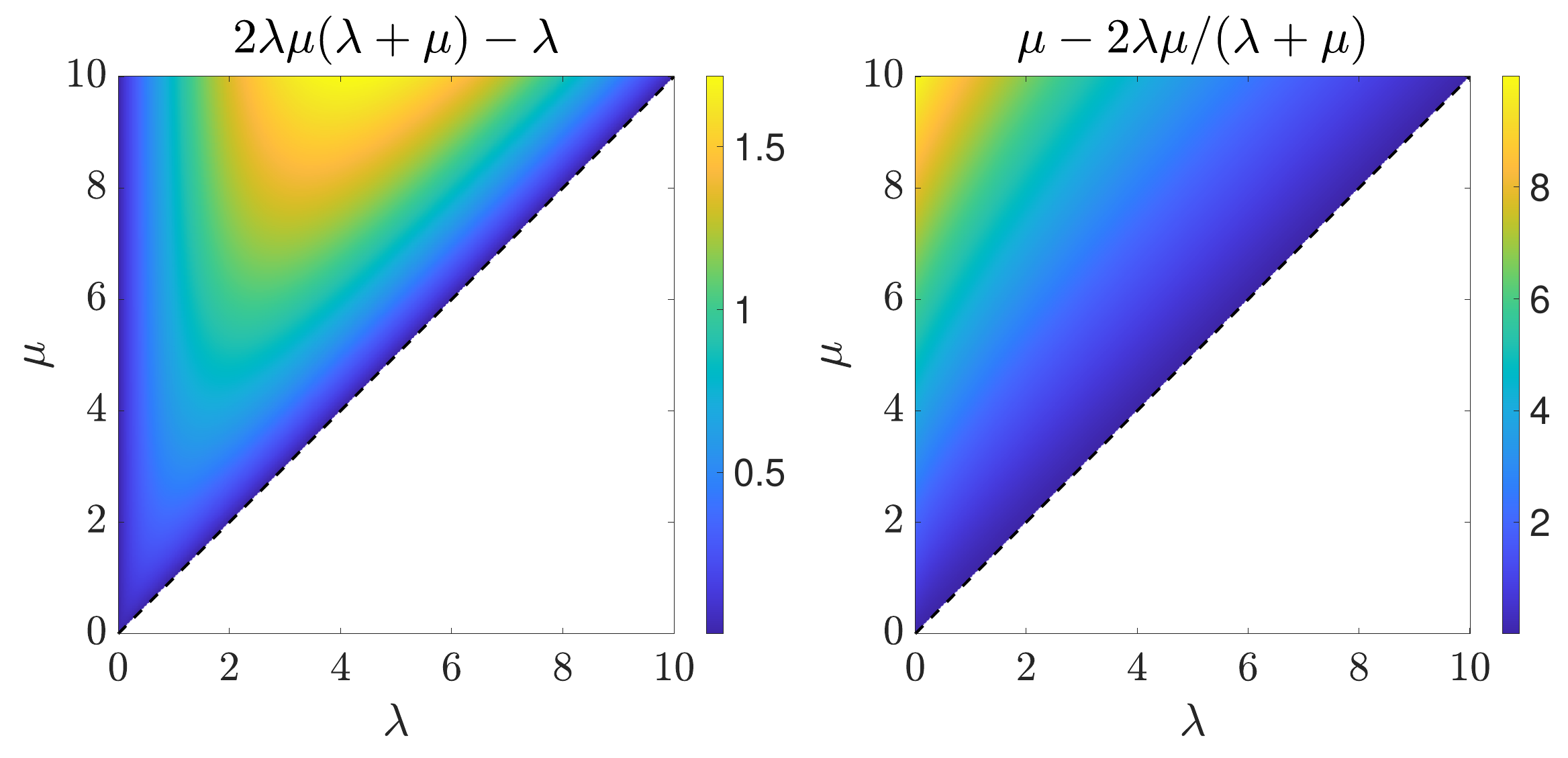}
\caption{\label{fig_bias}\textbf{Binary case.} Asymptotic bias of the classical MLEs \( \hat{\lambda}_t \) (left) and \( \hat{\mu}_t \) (right), as functions of pairs $(\lambda,\mu)$ with $\lambda<\mu$ (subcritical case).}
\end{figure}

We start with the standard subcritical birth-and-death process with $p_2=1$ (hence $m=2$ and $\sigma^2=0$).  
Recall from Remark~\ref{rem_int} that the limits $\lambda^\uparrow$ and $\mu^\uparrow$ of the classical ($Q$-consistent) MLEs \( \hat{\lambda}_t \) and \( \hat{\mu}_t \) satisfy
\[
\lambda \;<\; \lambda^\uparrow = \mu^\uparrow = \frac{2\lambda\mu}{\lambda+\mu} \;<\; \mu .
\]
In Figure~\ref{fig_bias} we show the asymptotic bias of the classical MLEs for $\lambda$ (left) and $\mu$ (right), as functions of pairs $(\lambda,\mu)$ with $\lambda<\mu$ ($\rho <0$).

We see that the bias is systematically larger for the estimator of $\mu$ than for that of $\lambda$.

The \emph{relative} asymptotic error induced by the $Q$-consistent estimators depends only on the ratio $\mu/\lambda$, and is given by the monotone increasing function 
\[ 
\dfrac{\lambda^\uparrow-\lambda}{\lambda}=\dfrac{\mu-\mu^\uparrow}{\mu}=\frac{(\mu/\lambda) - 1}{(\mu/\lambda) + 1},
\]
which tends to 1 as $\mu/\lambda\to\infty$. 
	
These observations highlight that the use of classical MLEs in subcritical birth-and-death processes can lead to substantial and systematic bias in the estimated birth and death rates, particularly when the ratio $\mu/\lambda$ is large.

\subsubsection*{Multiple-birth case}

We first make some observations on the relative asymptotic errors induced by the $Q$-consistent MLEs $\hat{\lambda}_t$, $\hat{\mu}_t$, and $\hat{p}_{k,t}$ ($k\geq 2$):
\[
\frac{\lambda^\uparrow-\lambda}{\lambda}=\frac{m-1}{\pi^\uparrow}, 
\qquad 
\frac{\mu-\mu^\uparrow}{\mu}=\frac{1}{\pi^\uparrow},
\]
and
\[
k\leq m:\quad 
\frac{p_k-p_k^\uparrow}{p_k}=\frac{m-k}{\pi^\uparrow+m-1}, 
\qquad 
k>m:\quad 
\frac{p_k^\uparrow-p_k}{p_k}=\frac{k-m}{\pi^\uparrow+m-1}.
\]
From Lemma \ref{lem:pi-up}, we can write
\[
\pi^\uparrow = 1+\frac{\sigma^2+m(m-1)}{(\mu/\lambda)-(m-1)}, 
\qquad 
\pi^\uparrow+m-1 = m+\frac{\sigma^2+m(m-1)}{(\mu/\lambda)-(m-1)},
\]
with $(\mu/\lambda)-(m-1)>0$ in the subcritical case. 

It follows that $(\pi^\uparrow)^{-1}$ and $(\pi^\uparrow+m-1)^{-1}$ are decreasing in $m$ and $\sigma^2$, and increasing in $\mu/\lambda$. As a consequence, all relative asymptotic errors increase with $\mu/\lambda$ and decrease with $\sigma^2$. In addition, the errors for $\hat{\mu}_t$ and for $\hat{p}_{k,t}$ with $k>m$ also decrease with $m$. The dependence on $m$ is less straightforward for the errors of $\hat{\lambda}_t$ and of $\hat{p}_{k,t}$ with $k\leq m$: since in the subcritical case $2\leq m<(\mu/\lambda)+1$, these errors may either be strictly decreasing in $m$ over this interval, or show a unimodal behaviour (increasing then decreasing), depending on the combination of $\mu/\lambda$ and $\sigma^2$.

Next, we consider a concrete example of a subcritical linear birth-and-death process with parameters  $$\lambda=2,\qquad \mu=5,\qquad (p_2,p_3,p_4)=(0.6,0.1,0.3),$$ which yields the offspring mean $m=2.7$ and variance $\sigma^2=0.81$. This process has mean growth rate $\rho=-1.6$, and the corresponding $Q$-process has asymptotic mean population size $\pi^\uparrow=7.75$. For this model, the limits of the $Q$-consistent MLEs are
$$\lambda^\uparrow=2.4387,\qquad \mu^\uparrow=4.3548,\qquad (p_2^\uparrow,p_3^\uparrow,p_4^\uparrow)=(0.5556,0.1032,0.3413),$$
with relative asymptotic errors 
\[
\frac{\lambda^\uparrow-\lambda}{\lambda}=0.219, \quad 
\frac{\mu-\mu^\uparrow}{\mu}=0.129, \quad
\frac{p_2-p_2^\uparrow}{p_2}=0.074, \quad
\frac{p_3^\uparrow-p_3}{p_3}=0.032, \quad
\frac{p_4^\uparrow-p_4}{p_4}=0.138.
\]
This illustrates the non-negligible bias of the classical estimators.

\begin{figure}[h!]
\centering\includegraphics[width=1\textwidth]{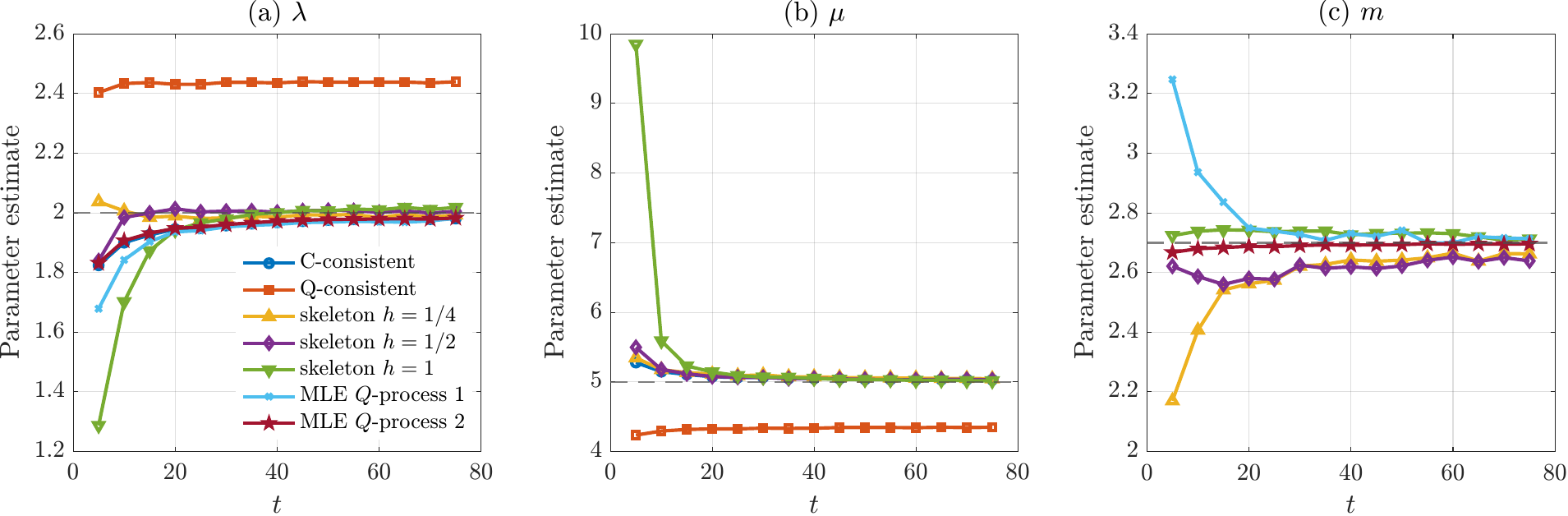}
\caption{\textbf{Multiple-birth case ($Z_0=5$).} Median estimates of the parameters $\lambda$, $\mu$, and $m$ as functions of the observation time $t$, based on 1500 simulated trajectories with initial population size $Z_0=5$. Each curve corresponds to a different estimator: $C$-consistent ($\tilde{\lambda}_t,\tilde{\mu}_t$), $Q$-consistent ($\hat{\lambda}_t,\hat{\mu}_t$), skeleton $h$ (the $\delta$-skeleton approach with step size $h$), and MLE $Q$-process (version~1: $m$ estimated via fixed point; version~2: finite offspring support, $m$ estimated via the $p_k$'s). Grey dashed lines indicate the true parameter values.
}
  \label{fig:consistency-params}
\end{figure}
\begin{figure}[h!]
\centering\includegraphics[width=1\textwidth]{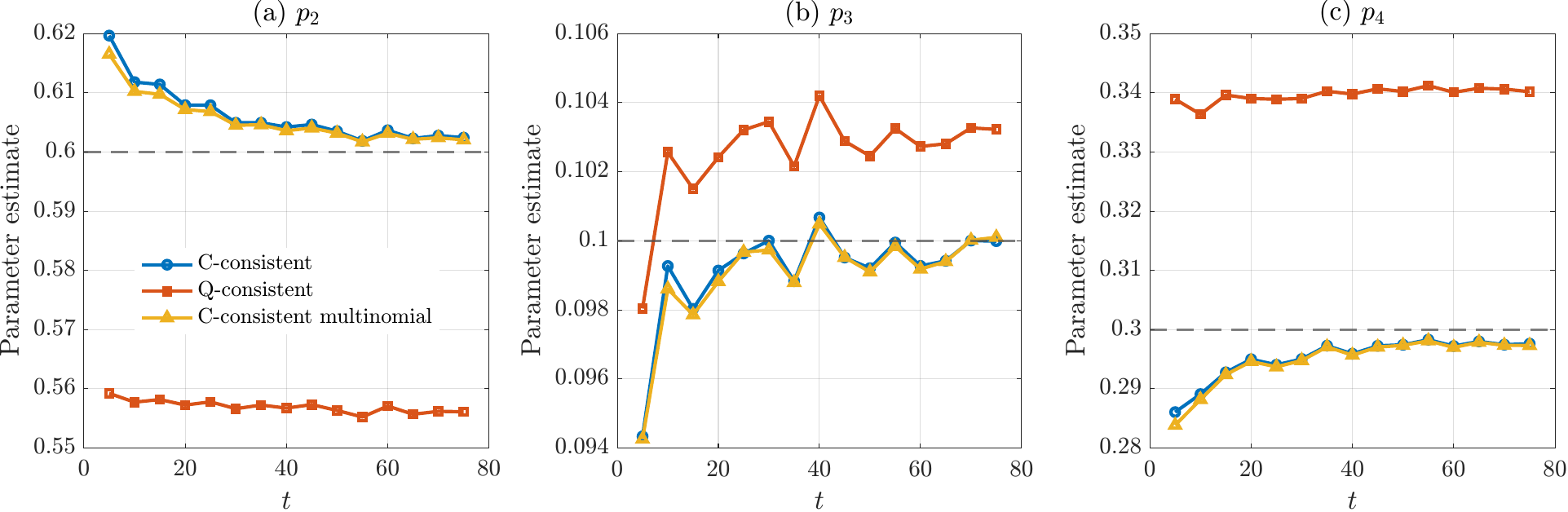}
\caption{\textbf{Multiple-birth case ($Z_0=5$).} Median estimates of the offspring probabilities $p_2,p_3,p_4$ as functions of the observation time $t$, based on 1500 simulated trajectories with $Z_0=5$. Three variants are compared: $C$-consistent $\tilde{p}_{k,t}$, $Q$-consistent $\hat{p}_{k,t}$, and the $C$-consistent normalised `multinomial' version $\bar{p}_{k,t}$. Grey dashed lines indicate the true parameter values.}
  \label{fig:consistency-pk}
\end{figure}

\begin{figure}[h!]
\centering\includegraphics[width=1\textwidth]{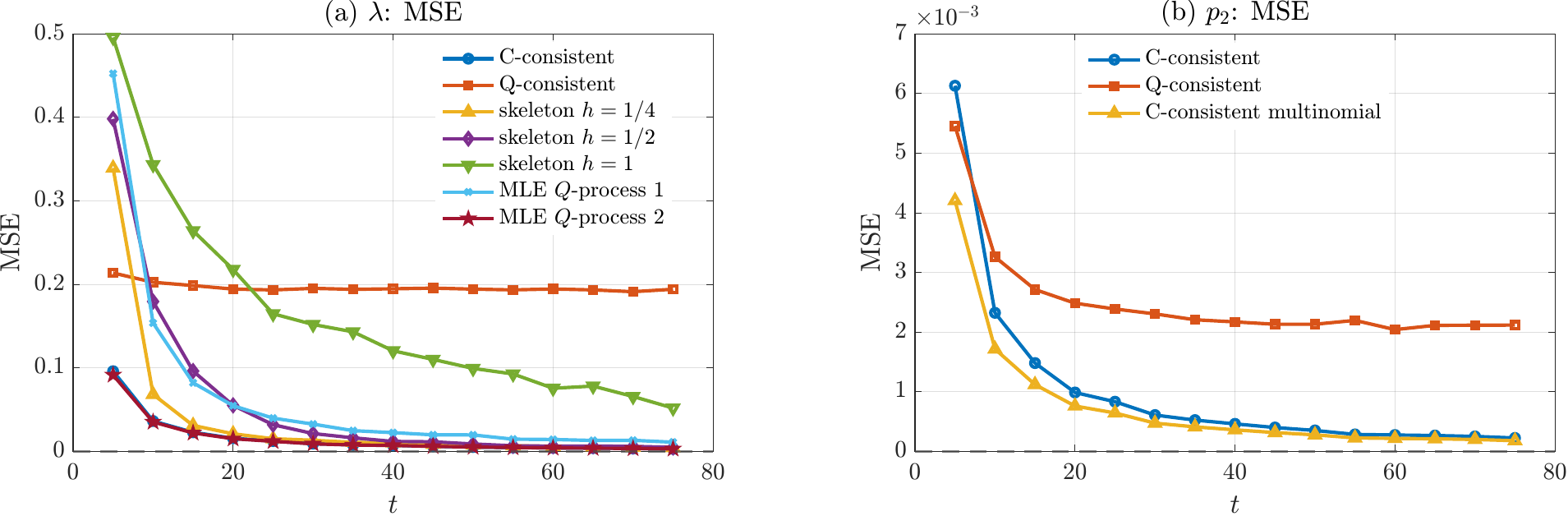}
\caption{\textbf{Multiple-birth case ($Z_0=5$).} Mean squared errors (MSE) of the estimators of (a)~$\lambda$ and (b)~$p_2$ as functions of $t$.}
  \label{fig:mse}
\end{figure}

Figures~\ref{fig:consistency-params} and~\ref{fig:consistency-pk} illustrate the convergence of the estimators introduced in Sections~\ref{sec:C-con}--\ref{rem-est-m} for $\lambda$ (Figure~\ref{fig:consistency-params}~(a)), $\mu$ (Figure~\ref{fig:consistency-params}~(b)), $m$ (Figure~\ref{fig:consistency-params}~(c)), and $(p_2,p_3,p_4)$ (Figure~\ref{fig:consistency-pk}), based on 1500 simulated trajectories with initial population size $Z_0=5$, from length $t=5$ to $t=75$. Figure~\ref{fig:mse} reports the mean squared errors (MSE) for selected estimators, namely those for $\lambda$ and $p_2$. To generate long non-extinct trajectories of subcritical birth-and-death processes we used a multilevel splitting method, similar to that used in \cite{BHM22} but adapted to the continuous-time setting.

The plots confirm the convergence of all $C$-consistent estimators to the true parameter values as $t$ increases. The estimators $\tilde{\lambda}_t$ and ${\bar{\lambda}}_t$, although not theoretically equivalent, yield almost indistinguishable estimates in practice. The two variants of $C$-consistent estimators for $(p_2,p_3,p_4)$ also produce very close results, with the multinomial version converging slightly faster. The $\delta$-skeleton approach improves in terms of MSE as the step size $h$ decreases. Among the $C$-consistent estimators derived via the $Q$-process MLE approach, the fixed-point version is generally less efficient than the finite-support version where $m$ is estimated via the $p_k$s.

\begin{figure}[ht]
 \centering\includegraphics[width=1\textwidth]{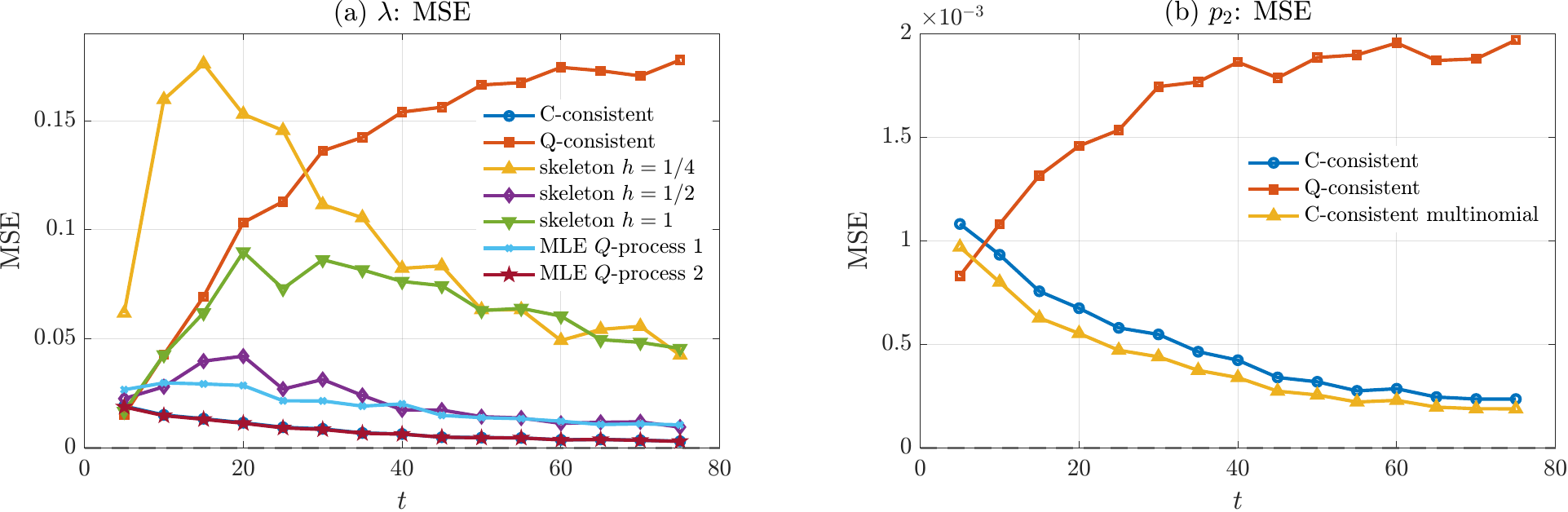}
 \caption{\textbf{Multiple-birth case ($Z_0=200$).} Mean squared errors (MSE) of the estimators of (a)~$\lambda$ and (b)~$p_2$ as functions of $t$. \label{fig:appc}}
 \end{figure}

Finally, Figure \ref{fig:appc} reports MSE results for processes starting with $Z_0=200$. In that case, the $Q$-consistent estimators have smaller MSE than certain $C$-consistent estimators for very short observation windows (i.e., when the population size is still far from extinction), whereas the $C$-consistent estimators outperform them as the observation windows increase.

\section{Proofs}\label{sec:proofs}
 
To prove our results, for each $t\geq 0$ we place ourselves in the probability space of a MEXIT coupling of the non-homogeneous process $(Z_s^{(t)})_{0\leq s \le t}$ and the $Q$-process $(Z_s^\uparrow)_{0\leq s \le t}$ as defined in Appendix~\ref{sec:coupling}, with probabilities $(\widehat{\mbP}^{(t, \uparrow)})_{t \ge 0}$ and expectations $(\widehat{\mbE}^{(t, \uparrow)})_{t \ge 0}$, which we simplify to $\mathbb{P}$ and $\mathbb{E}$ for convenience.
Recall that $\zeta_t$ is the uncoupling time defined in \eqref{uncoupling-time}, and $C(i,q)$ is the constant defined in Proposition \ref{thm:mexit}\emph{(iii)}.

For $t \ge 0$ and $q>0$, we define the events
\begin{align*}
A_{1, t} &:= \{\zeta_t > t - C(i, q)\log t\}\\
A_{2, q, t} &:= \{Z_s^\uparrow \leq t^q \text{ for all } s \in [t - C(i, q)\log t, t]\}\\
A_{3, q, t} &:= \{Z_s^{(t)} \leq t^q \text{ for all } s \in [t - C(i, q)\log t, t]\}.
\end{align*}
The proofs of the results in Sections \ref{sec:C-con} and \ref{sec:pk} rely on a series of lemmas which we state below.


\begin{lemma}\label{lem-events}
For all $i>0$ and $q>0$,  $\mbP_i(A_{1, t}) \to 1$, $\mbP_i(A_{2, q, t}) \to 1$ and $\mbP_i(A_{3, q, t}) \to 1$ as $t \to \infty$.
\end{lemma}

\noindent \textbf{Proof.}
From Proposition \ref{thm:mexit} (iii), we have $\mbP_i(A_{1, t}) \to 1$.
To show $\mathbb{P}_i	(A_{2,q,t}) \to 1$, we show that $\mathbb{P}_i	(A_{2,q,t}^c) \to 0$. Using Markov's inequality and the Markov property, we have
\begin{align}
\mathbb{P}_i \left( \sup_{s \in [t - C(i,q) \log t, t]} Z^\uparrow_s > t^q \right)
&\leq t^{-q} \, \mathbb{E}_i \left[ \sup_{s \in [t - C(i,q) \log t, t]} Z^\uparrow_s \right] \nonumber \\
&= t^{-q} \, \mathbb{E}_i \left[ \mathbb{E} \left( \left. \sup_{0\leq s \leq C(i,q) \log t} Z^\uparrow_s \right| Z^\uparrow_{t - C(i,q) \log t} \right) \right].
\label{eq:markov_bound}
\end{align}

By the spine decomposition described in Section \ref{sec:cond}, we know that, if we start with $i$ individuals, one uniformly chosen individual initiates a copy of the spine and the remaining $i - 1$ individuals initiate copies of the original process. Thus,
\begin{multline}
\mathbb{E} \left[ \left. \sup_{0\leq s \leq C(i,q) \log t} Z^\uparrow_s \right| Z^\uparrow_{t - C(i,q) \log t} \right] \\
\leq
\mathbb{E} \left[ \sup_{0\leq s \leq C(i, q) \log t} Z^\uparrow_s \right]
+ \left( Z^\uparrow_{t - C(i,q) \log t} - 1 \right) \mathbb{E} \left[ \sup_{0\leq s \leq C(i,q) \log t} Z_s \right].
\label{eq:spine_split}
\end{multline}

Let us consider the first of the two expectations on the right-hand side above. Again using the spine decomposition \eqref{spine_dec}, we have
\[
Z^\uparrow_s = 1 + \sum_{\ell = 1}^{N_s}  \sum_{\substack{j=1\\ j \neq i^*}}^{\tilde{\xi}_\ell} Z^{(j)}_{s - T_\ell},
\]
where $N_s\sim \textrm{Poi}(m\lambda s)$ denotes the number of birth events along the spine until time $s$, $T_1,\ldots, T_{N_s}$ are the times of these birth events, $\tilde{\xi}\sim (\tilde{p}_k)_{k\geq 2}$ is the size biased offspring distribution, and $i^* \sim {\rm Unif}\{1, \dots, \tilde\xi_\ell\}$ denotes the index of the spine particle.
Since $(N_s)_{s \geq 0}$ is non-decreasing and the summands are non-negative, we get
\[
\sup_{0\leq s \leq C(i,q) \log t} Z^\uparrow_s \leq 
1 + \sum_{\ell = 1}^{N_{C(i,q) \log t}}  \sum_{\substack{j=1\\ j \neq i^*}}^{\tilde{\xi}_\ell} \sup_{0\leq s \leq C(i,q) \log t} Z^{(j)}_{s - T_\ell}.
\]
Taking expectations and following the same argument as in the proof of Lemma \ref{lem:pi-up}, we obtain
\begin{equation}
\mathbb{E} \left[ \sup_{0\leq s \leq C(i,q) \log t} Z^\uparrow_s \right] \leq 1 + \lambda (\sigma^2 + m(m-1)) \int_0^{C(i,q) \log t} \mathbb{E} \left[ \sup_{0\leq s \leq C(i, q) \log t} Z_{s - u} \right] du.
\label{eq:spine_bound}
\end{equation}

Since the process is subcritical and the offspring distribution has finite variance, it follows that
\[
M := \mathbb{E} \left[ \sup_{s \geq 0} Z_s \right] < \infty;
\]see for example \cite{janson1986moments}.
Using this in \eqref{eq:spine_bound}, we get
\begin{equation}
\mathbb{E} \left[ \sup_{s \leq C(i, q) \log t} Z^\uparrow_s \right] \leq 1 + \lambda (\sigma^2 + m(m-1))\,M\, C(i,q) \log t.
\label{eq:spine_final_bound}
\end{equation}
Combining \eqref{eq:spine_split} and \eqref{eq:spine_final_bound}, taking expectations, and using \eqref{ab} and \eqref{M1-up}, we obtain
\begin{eqnarray*}\lefteqn{
\mathbb{E}_i \left[ \mathbb{E} \left( \left. \sup_{0\leq s \leq C(i,q) \log t} Z^\uparrow_s \right| Z^\uparrow_{t - C(i,q) \log t} \right) \right]}\\
&\leq & 1 + \lambda (\sigma^2 + m(m-1)) M\, C(i,q) \log t + \mathbb{E}_i \left[ Z^\uparrow_{t - C(i,q) \log t} - 1 \right] M \\
&=& 1 + \lambda (\sigma^2 + m(m-1)) M\, C(i,q) \log t \\&&+ \lambda (\sigma^2 + m(m-1)) \frac{e^{\rho(t - C(i,q) \log t)} - 1}{\rho}\,M+(i-1) e^{\rho(t - C(i,q) \log t)}\, M.
\end{eqnarray*}

Returning to \eqref{eq:markov_bound} and dividing the right-hand side by $t^q$, then letting $t \to \infty$, yields the result:
\[
\mathbb{P}_i \left( \sup_{s \in [t - C(i,q) \log t, t]} Z^\uparrow_s > t^q \right) \to 0,
\]
which is what was required.

Finally, to show $\mbP_i(A_{3, q, t}) \to 1$, we observe that, thanks to stochastic domination (see for instance \cite[Corollary 4]{BHM22}), we have $\mbE_i[Z_s^{(t)}] \le \mbE_i[Z_s^\uparrow]$ for all $s \le t$, $t \ge 0$ and $i \ge 1$. In particular, we have 
\[\mathbb{E}_i \left[ \sup_{s \in [t - C(i,q) \log t, t]} Z^{(t)}_s \right]\leq \mathbb{E}_i \left[ \sup_{s \in [t - C(i,q) \log t, t]} Z^\uparrow_s \right],\]and we then apply
the same arguments as above, which concludes the proof. \hfill $\square$

\medskip

\begin{lemma}\label{lem:uparrow_conv}For any initial population size, as $t\to\infty$,
\begin{itemize}
\item[(i)]  ${\tau_t^\uparrow}/{t}\xrightarrow{p}\pi^\uparrow$;
\item[(ii)] ${b_t^\uparrow}/{t}\xrightarrow{p}\lambda(\pi^\uparrow - 1) + \lambda m$;
\item[(iii)] ${d_t^\uparrow}/{t}\xrightarrow{p}\mu(\pi^\uparrow - 1)$.
\end{itemize}
\end{lemma}

\medskip

\noindent \textbf{Proof.}
 \emph{(i)} The convergence of \( \tau_t^\uparrow / t \) to \( \pi^\uparrow \) follows directly from the ergodic theorem for the positive recurrent Markov chain \( Z^\uparrow \) which states that the time average \( \tfrac1t\int_0^t Z_u^\uparrow\,du \) converges almost surely to the mean of the stationary distribution of \( Z^\uparrow \); see, for example, \cite[Section~5.5]{resnick2013adventures}.

\smallskip

\noindent \emph{(ii)} To analyse \( b_t^\uparrow / t \), we apply the spine decomposition. In the case where the initial population size $i=1$, this allows us to write
\begin{equation}\label{spdec2}
\dfrac{b_t^\uparrow}{t} = \frac{1}{t}\sum_{i=1}^{N_t} \left(1+\sum_{j=1}^{\tilde{\xi}_i - 1} b^{(i, j)}\right) - \frac{1}{t}\sum_{\ell=1}^{Z_t^\uparrow-1}\tilde b_{t}^{(\ell)},
\end{equation}
where \( N_t \) is the number of birth events along the spine up to time \( t \),
 \( \tilde{\xi}_i \) is the number\footnote{Note that we have slightly abused notation here by assuming that the first $\tilde\xi_i - 1$ are the non-spine particles however, since they are all i.i.d. and the spine is chosen uniformly, it makes no difference to the subsequent analysis.} of offspring generated at the \( i \)-th birth event along the spine,
   \( b^{(i, j)} \) is the total number of birth events until extinction in the \( j \)-th subcritical process born in the $i$-th birth event along the spine, and $\tilde b_{t}^{(\ell)}$ denotes the total number of birth events until extinction in the process initiated from $\ell$-th non-spine particle alive at time $t$.

Define \( Y_i := 1 + \sum_{j=1}^{\tilde{\xi}_i - 1} b^{(i, j)} \), for \( i = 1, \dots, N_t \). Since the \( Y_i\) are i.i.d. and also independent of $\tilde\xi_i$, $\mathbb{E}[\tilde{\xi} - 1] = \frac{\sigma^2}{m} + m - 1$ and $\mathbb{E}[b] = -\frac{\lambda}{\rho}$, it follows that
\[
\mathbb{E}[Y] = 1 + \mathbb{E}[\tilde{\xi} - 1] \cdot \mathbb{E}[b] = 1 - (\tfrac{\sigma^2}{m} + m - 1)\frac{\lambda}{\rho}.
\]
 Since \( N_t \to \infty \) almost surely as \( t \to \infty \), the law of large numbers gives
\[
\frac{1}{N_t} \sum_{i=1}^{N_t} Y_i \xrightarrow{\text{a.s.}} \mathbb{E}[Y], \quad t \to \infty.
\]
Similarly,
$
\frac{N_t}{t} \xrightarrow{a.s.} \lambda m,
$ as $t \to \infty$. Hence, by the Continuous Mapping Theorem and Lemma~\ref{lem:pi-up}, we obtain
\[
\frac{1}{t} \sum_{i=1}^{N_t} Y_i = \frac{N_t}{t} \cdot \frac{1}{N_t} \sum_{i=1}^{N_t} Y_i \xrightarrow{a.s.} \lambda m \cdot \mathbb{E}[Y] = \lambda m + \lambda(\pi^\uparrow - 1),
\]
as \( t \to \infty \).

For the second term in \eqref{spdec2}, using \eqref{ab} and \eqref{M1-up} we see that \( \mbE[Z_t^\uparrow] \) is  uniformly bounded in $t$, and since \( \mathbb{E}[\tilde b_{t}^{(\ell)}] < \infty \), it follows that
\[
\frac{1}{t} \sum_{\ell=1}^{Z_t^\uparrow-1} \tilde b_{t}^{(\ell)} \xrightarrow{p} 0,
\]
since it converges to zero in mean.

Finally, if the initial population size is $i>1$, it suffices to add to \eqref{spdec2} a term which is bounded above by $\frac1t \sum_{\ell=1}^{i-1} b^{(\ell)} \xrightarrow{p} 0$ as $t \to \infty$, where, here, the $b^{(\ell)}$ are i.i.d. copies of the number of birth events in a subcritical birth-death process until extinction.

\smallskip

\noindent \emph{(iii)} The arguments for the proof of ${d_t^\uparrow}/{t}\xrightarrow{p}\mu(\pi^\uparrow - 1)$ are similar to those used in \emph{(ii)}, noting that
$$\dfrac{d_t^\uparrow}{t} = \frac{1}{t}\sum_{i=1}^{N_t} \left(\sum_{j = 1}^{\tilde{\xi}_i - 1} d^{(i, j)}\right) - \frac{1}{t}\sum_{\ell=1}^{Z_t^\uparrow-1} \tilde d_{t}^{(\ell)},$$ 
where the terms in this expression are defined analogously to those above for the number of death events. 
\hfill $\square$

\bigskip

We use the notation $\theta_t^{(t)}$ to denote an estimator or statistics $\theta_t$ based on observations of the inhomogeneous process $(Z_u^{(t)})_{0\leq u\leq t}$, that is, the process $Z$ conditional on $\{Z_t > 0\}$. 

\begin{lemma}\label{lem:diff_t_uparrow}For any initial population size and $p>0$, as $t\to\infty$,
\begin{itemize}
\item[(i)]  $\left(\tau_t^{(t)} - \tau_t^\uparrow\right)/t^p\xrightarrow{p} 0;$
\item[(ii)] $\left(b_t^{(t)} - b_t^\uparrow\right)/t^p\xrightarrow{p} 0$;
\item[(iii)] $\left(d_t^{(t)} - d_t^\uparrow\right)/t^p\xrightarrow{p} 0$.
\end{itemize}
\end{lemma}

\noindent \textbf{Proof.}
 By Lemma \ref{lem-events}, it suffices to consider the sequences of random variables on $A_{1, t},$ 
$A_{2, q, t},$ and $
A_{3, q, t}$ for well-chosen $q$.

\smallskip

\noindent \emph{(i)} On $A_{1, t},$ 
$A_{2, q, t},$ and $
A_{3, q, t}$ with $0<q<p$, for any initial population size $i > 0$, we have 
\begin{align*}
  \frac{1}{t^p}\mbE_i\big[|\tau_t^{(t)} - \tau_t^\uparrow|\big] 
  &= \frac{1}{t^p}\mbE_i\left[\left|\int_0^t Z_s^{(t)} \dd s - \int_0^t Z_s^\uparrow \dd s \right|\right]\\
  & = \frac{1}{t^p}\mbE_i \left[\left|\int_{\zeta_t}^t \big(Z_s^{(t)} -  Z_s^\uparrow \big)\, \dd s \right|\right] \\
  & \le \frac{1}{t^p} \int_{t-C(i, q)\log t}^t \mbE_i\big[|Z_s^{(t)} -  Z_s^\uparrow|\big]\dd s \\
  & \le \frac{2C(i, q) \log t}{t^{p-q}},
\end{align*}
which converges to $0$ as $t \to \infty$.

\smallskip

\noindent \emph{(ii)} We will show that, on the events $A_{1, t},$ 
$A_{2, q, t},$ and $
A_{3, q, t}$ with $0<q<p$, we have
\begin{equation}\label{ine2}
\mbE_i\left[\left|b_t^{(t)}-b_t^{\uparrow}\right|\right]
= O(t^q\log t), \quad t \to \infty.
\end{equation}
Since convergence in $L^1$ implies convergence in probability, the result follows.

To this end, first note that 
\[
\mbE_i\left[\left|b_t^{(t)}-b_t^{\uparrow}\right|\right]
\leq \mbE_i\left[B^{(t)}_{[t-C(i,q)\log t,t]}\right]+\mbE_i\left[B^\uparrow_{[t-C(i,q)\log t,t]}\right],
\]
where $B^{(t)/\uparrow}_{[t-C\log t,t]}$ denotes the number of birth events in the process $Z^{(t)/\uparrow}$ during the time period $[t-C\log t,t]$ and we have written (and will do so throughout the rest of the proof) $C = C(i, q)$ for brevity. This follows from the fact that, on $A_{1, t},$ the processes $(Z_s^{(t)})_{0\leq s\leq t-C\log t}$ and $(Z_s^{\uparrow})_{0\leq s\leq t-C\log t}$ are equal. 

Next we show that $\mbE_i[B^\uparrow_{[t-C\log t,t]}] = O(t^q \log t)$ as $t \to \infty$. First note that, on the event $A_{2, q, t}$, the number of branches is bounded above by $t^q$. Moreover, the birth rate is bounded above by $\lambda m$. Hence, the number of birth events in the time interval $[t-C\log t,t]$ is dominated by a Poisson random variable with  mean $\lambda m\,t^{q}\,C\,\log t$, and hence the claim follows. 

Finally, we show that this also holds for the conditioned process. This is more involved since, in this case, the birth rate is time-dependent and cannot be so easily bounded from above.

In general, the time-inhomogeneous transition rates at time $s$ of a continuous-time Markov chain $Z$ conditioned on non-absorption at time $t\geq s$ are given by
$$Q^*_{i,j}(s;t)=Q_{ij}\,\dfrac{h(j,s,t)}{h(i,s,t)}, \qquad i,j\neq 0,$$
where $h(i,s,t):=\mbP(Z_t> 0\mid Z_s=i)$. In our case, 
$$h(i,s,t)=1-\mbP(Z_t =0\mid Z_s=i)=1- F(t-s)^i,$$
where $F(\tau):=\mbP(Z_\tau=0\mid Z_0=1)$. 

So the total birth rate at time $s$ and at population size $z$ in the time-inhomogeneous process $(Z_s^{(t)})_{t-C\log t\leq s\leq t}$ is given by
\begin{equation}\lambda^*_z(s; t)=\sum_{k\geq 2} \lambda\, p_k \,z\, \dfrac{\left[1-F(t-s)^{z-1+k}\right]}{1-F(t-s)^z}=\lambda \,z\,\dfrac{1-F(t-s)^{z-1}\, P(F(t-s))}{1-F(t-s)^z},
\end{equation} 
where $P(s):=\sum_{k\geq 2} p_k s^k$ is the p.g.f.\ of the offspring distribution.
We make the following observations, assuming $t$ is fixed:\begin{itemize}
\item The function $F(t-s)$ is a decreasing function of $s$ over $[t-C\log t, t]$ with $F(0)=0$.
\item The rate $\lambda^*_z(s; t)$ is also a decreasing function of $s$ over $[t-C\log t, t]$ with $\lambda^*_z(t;t)= \lambda\, z$ (the total birth rate in the original unconditional process). We thus have, for all $z\geq 1$ and $s\in [t-C\log t, t]$,
$$ \lambda^*_z(s;t)\leq \lambda^*_z(t-C\log t;t)
=\lambda \,z\,\dfrac{1-F(C\log t)^{z-1}\, P(F(C\log t))}{1-F(C\log t)^z}.$$
\item The rate $\lambda^*_z(s;t)$ is an increasing function of $z$. Therefore, on $A_{3, q, t}$, for all $z\geq 1$ and $s\in [t-C\log t, t]$, the total birth rate is bounded above by
$$ \lambda^*_z(s;t)\leq  \lambda \,t^q\,\dfrac{1-F(C\log t)^{t^q-1}\, P(F(C\log t))}{1-F(C\log t)^{t^q}}=: U(t,q).$$
So the number of birth events over $[t-C\log t,t]$ is dominated by a Poisson random variable with mean $U(t,q)\, C\log t$.
\end{itemize}  
 We now study the asymptotic behaviour of $U(t,q)$ as $t \to \infty$. In this asymptotic regime, it is well known that $$F(C\log t)\sim 1-Ke^{\rho\, C\log t}=1-Kt^{\rho\,C}, \quad t \to \infty$$
 for some constant $K > 0$. Therefore   
\begin{equation}\label{eq:Uasymp}
U(t,q)\sim  \lambda \,t^{q}\,\frac{1 - \left(1 - Kt^{\rho C}\right)^{t^q  - 1}\,P(1- Kt^{\rho C})}{1 - \left(1 - Kt^{\rho C}\right)^{t^q}}, \quad t \to \infty.
\end{equation}
Now, for $t$ sufficiently large (so that $Kt^{\rho C} < \frac12$), the Taylor expansion of $\log(1-x)$ yields 
\[
 t^q\log(1-Kt^{\rho C}) = -Kt^{q + \rho C} + t^qO(t^{2\rho C}).
\]
Choosing $0 < q < \min\{|\rho|C, p\}$, the above entails that
\[
t^q\log(1-Kt^{\rho C}) = -Kt^{q + \rho C} + o(t^{q + \rho C}). 
\]
Hence, 
\[
 (1-Kt^{\rho C})^{t^q} = 1 - Kt^{q + \rho C} + o(t^{q + \rho C}), \quad t \to \infty,
\]
which follows from the fact that ${\rm e}^y = 1 + y + o(y)$ as $y \to 0$. Similarly, we also have
\[
 (1-Kt^{\rho C})^{t^q - 1} = 1 - K(t^{q}-1)t^{\rho C} + o((t^{q}-1)t^{\rho C}), \quad t \to \infty.
\]
Substituting this back into \eqref{eq:Uasymp} yields
\[
  U(t, q) \sim \lambda \,t^{q}\frac{K(t^q - 1)t^{\rho C} + o((t^q - 1)t^{\rho C})}{Kt^{q + \rho C} + o(t^{q + \rho C})}, \quad t \to \infty.
\]
Dividing both the numerator and denominator by $Kt^{q + \rho C}$, the right-hand side is equal to
\[
  \lambda \,t^{q}\frac{1 - t^{-q} + o(t^{-q})}{1 + o(1)} = \lambda(t^{q} - 1) + o(1), \quad t \to \infty.
\]
This implies that $\mbE[B^{(t)}_{[t-C\log t,t]}] =  O(t^{q}\log t)$, as $t \to \infty$, which completes the proof.

\smallskip
\noindent \emph{(iii)} The arguments for the proof of $\left(d_t^{(t)} - d_t^\uparrow\right)/t^p\xrightarrow{p} 0$ are similar to those used in \emph{(ii)}, noting that the total death rate at time $s$ and at population size $z$ in the time-inhomogeneous process $(Z_s^{(t)})_{t-C\log t\leq s\leq t}$ is given by
\begin{equation}\mu^*_z(s;t)=\mu \,z\,\dfrac{1-F(t-s)^{z-1}}{1-F(t-s)^z},
\end{equation}which is an increasing function of $s$ over $[t-C\log t, t]$, bounded above by $\mu^*_z(t;t)=\mu\,z$ (the total death rate in the original unconditional process).
 \hfill $\square$

\begin{cor}\label{lem:tautt}
For any initial population size, as $t\to\infty$,
  $\tau_t^{(t)}/{t}\xrightarrow{p}\pi^\uparrow$.
\end{cor}

\noindent \textbf{Proof.} This is a direct consequence of the triangle inequality, Lemma \ref{lem:uparrow_conv}\emph{(i)} and Lemma \ref{lem:diff_t_uparrow}\emph{(i)} with $p=1$.\hfill $\square$

\subsection{Proof of the results in Section \ref{sec:C-con}}

\noindent\textbf{Proof of Proposition \ref{prop:Con_Qprocess}} (Consistency of $\hat{\lambda}_t^\uparrow$ and $\hat{\mu}^\uparrow_t$).
The proofs of \eqref{Q-con-lam2} and \eqref{Q-con-mu2} follow directly from the definition of the estimators $\hat\lambda_t^\uparrow$ and $\hat\mu_t^\uparrow$ in \eqref{MLEsQ}, Lemma \ref{lem:uparrow_conv}, the Continuous Mapping Theorem, and the definition of the limits $\lambda^\uparrow$ and $\mu^\uparrow$ in \eqref{lam_mu_uparrow}. \hfill $\square$
 
 \bigskip

\noindent \textbf{Proof of  Proposition \ref{prop:Qcon}} ($Q$-consistency of \( \hat{\lambda}_t \) and \( \hat{\mu}_t \)).
Since the proofs for \eqref{Q-con-lam} and \eqref{Q-con-mu} are almost identical, we focus on proving the result for $\hat\lambda_t$ only. 

\smallskip

Thanks to Proposition \ref{prop:Con_Qprocess}, it remains to show that $\hat\lambda_t^{(t)} - \hat\lambda_t^\uparrow \to 0$ in probability as $t \to \infty$. From the definitions of the estimators, we have
\begin{equation}
  \hat\lambda_t^{(t)} - \hat\lambda_t^\uparrow 
  = \frac{\frac{1}{t^2}(b_t^{(t)}\tau_t^\uparrow - b_t^\uparrow\tau_t^{(t)})}{\frac{1}{t^2}\tau_t^\uparrow \tau_t^{(t)}}.
  \label{eq:diff}
\end{equation}
We start by treating the denominator. We have
\begin{equation}\label{bb}
  \frac{1}{t^2}\tau_t^\uparrow \tau_t^{(t)}-(\pi^\uparrow)^2 = \frac{\tau_t^\uparrow}{t}\frac{\big(\tau_t^{(t)} - \tau_t^\uparrow\big)}{t} + \left(\frac{1}{t^2}(\tau_t^\uparrow)^2 - (\pi^\uparrow)^2\right).
\end{equation}
By Lemma \ref{lem:uparrow_conv}\emph{(i)}, we have that $\tau_t^\uparrow/t$ converges in probability to $\pi^\uparrow$. Therefore, by the Continuous Mapping Theorem,  the second term in the RHS of \eqref{bb} converges to $0$ as $t \to \infty$, and $\frac{\tau_t^\uparrow}{t}$ in the first term converges to $\pi^\uparrow$. By Lemma \ref{lem:diff_t_uparrow}\emph{(i)} with $p=1$, the remaining factor $\left(\tau_t^{(t)} - \tau_t^\uparrow\right)/t\xrightarrow{p} 0$.
 
We conclude that \eqref{bb} converges in probability to zero, that is, the denominator in \eqref{eq:diff} converges to $(\pi^\uparrow)^2$ as $t \to \infty$. 

\smallskip

We now turn to the numerator:
\begin{equation}\label{cc}
  \frac{1}{t^2} \big(b_t^{(t)}\tau_t^\uparrow - b_t^\uparrow\tau_t^{(t)}\big)
  = \frac{1}{t^2}\tau_t^\uparrow\big(b_t^{(t)} - b_t^\uparrow \big) 
  + \frac{1}{t^2}b_t^\uparrow\big(\tau_t^\uparrow - \tau_t^{(t)}\big).
\end{equation}
By Lemma \ref{lem:uparrow_conv}\emph{(i)} and \emph{(ii)}, ${\tau_t^\uparrow}/{t}\xrightarrow{p}\pi^\uparrow$ and ${b_t^\uparrow}/{t}\xrightarrow{p}\lambda(\pi^\uparrow - 1) + \lambda m$,  and by Lemma \ref{lem:diff_t_uparrow}\emph{(i)} and \emph{(ii)} with $p=1$,  $\left(\tau_t^{(t)} - \tau_t^\uparrow\right)/t\xrightarrow{p} 0,$ and  $\left(b_t^{(t)} - b_t^\uparrow\right)/t\xrightarrow{p} 0$. Using the Continuous Mapping Theorem, we conclude that \eqref{cc}, and therefore the numerator in \eqref{eq:diff}, converges in probability to zero as $t \to \infty$, which completes the proof. \hfill$\square$

\bigskip

\noindent \textbf{Proof of Theorem \ref{thm:CLT}} (Asymptotic normality of $\tilde{\lambda}_t$ and $\tilde{\mu}_t$). We prove the result for $\tilde{\lambda}_t$; the proof for $\tilde{\mu}_t$ follows similar arguments.

From the definition of $\tilde\lambda_t$ and some simple manipulation, we have
\[
  \frac{\sqrt{t}(\tilde\lambda_t - \lambda)}{\sqrt{\lambda/(\pi^\uparrow + m-1)}} = \frac{b_t - \lambda(\tau_t + (m-1)t)}{\sqrt{\lambda}\sqrt{t(\pi^\uparrow + m - 1})}\frac{\pi^\uparrow + m-1}{(\tau_t/t + m-1)}.
\]
Thanks to Corollary \ref{lem:tautt} and Slutsky's Theorem, it is sufficient to show that, as $t \to \infty$,
\begin{equation*}
\hat{Y}^{(t)}_t:=\dfrac{b_t^{(t)}-\lambda (\tau_t^{(t)} +(m-1)\,t)}{\sqrt{t(\pi^\uparrow+m-1)}} \xrightarrow{d} Z,
\end{equation*}
where $Z \sim \mathcal N(0, \lambda)$.

To this end, we also define
\begin{equation*} \hat{Y}^{\uparrow}_t:=\dfrac{b_t^{\uparrow}-\lambda (\tau_t^{\uparrow} +(m-1)\,t)}{\sqrt{t(\pi^\uparrow+m-1)}}.
\end{equation*}
We will first show that $\hat{Y}_t^{(t)} - \hat{Y}_t^\uparrow \xrightarrow{p} 0$ and then that $\hat{Y}^{\uparrow}_t\xrightarrow{d} Z$, as $t \to \infty$. Slutsky's Theorem then yields the result. 

For the first part, we have
$$\hat{Y}^{(t)}_t-\hat{Y}^{\uparrow}_t=\dfrac{(b_t^{(t)}-b_t^{\uparrow})+\lambda\,(\tau_t^{\uparrow}-\tau_t^{(t)})}{\sqrt{t(\pi^\uparrow+m-1)}},$$
which converges to $0$ in probability thanks to Lemma \ref{lem:diff_t_uparrow} with $p=1/2$, and the Continuous Mapping Theorem.

For the second part, we note that the $Q$-process, $Z^{\uparrow}$, is a birth-and-death process whose \textit{total} birth rate and death rate at population size $r$ take the respective forms $$\lambda_r = f(r)\lambda\quad\text{and}\quad \mu_r = g(r)\mu,\quad\textrm{ with }\quad f(r) = (r - 1) + m \quad\text{and}\quad g(r) = r - 1.$$
By \cite[Equations (2) and (3)]{reynolds}, we can conclude that the counterparts of $\tilde{\lambda}_t$ and $\tilde{\mu}_t$ in the $Q$-process, defined by
\begin{equation}\label{tilde_up}
  \tilde\lambda_t^\uparrow := \frac{b_t^\uparrow}{\tau_t^\uparrow + (m-1)t}, \qquad 
  \tilde\mu_t^\uparrow := \frac{d_t^\uparrow}{\tau_t^\uparrow - t},
\end{equation}
are the MLEs for $\lambda$ and $\mu$, respectively (note that allowing for multiple births in this setting does not change the form of the likelihood for a continuously observed trajectory). Then, by \cite[Equation (5)]{reynolds}, we obtain
$$\sqrt{t(\pi^\uparrow+m-1)}\left(\frac{b_t^\uparrow}{\tau_t^\uparrow + (m-1)t}-\lambda\right) \xrightarrow{d} Z, \quad t \to \infty.$$ 
Finally, we note that 
\[
\hat{Y}_t^\uparrow = \sqrt{t(\pi^\uparrow+m-1)}\left(\frac{b_t^\uparrow}{\tau_t^\uparrow + (m-1)t}-\lambda\right) \frac{\tau_t^\uparrow/t + m-1}{\pi^\uparrow + m-1},
\]
which also converges in distribution to $Z$ thanks to the above combined with Lemma \ref{lem:uparrow_conv} and Slutsky's Theorem. 

Independence of the asymptotic normal distributions for the estimators for $\lambda$ and $\mu$ arises from the fact that the information matrix of the MLEs for $\lambda$ and $\mu$ is diagonal (see \cite{reynolds}).
\hfill$\square$

\bigskip
\noindent \textbf{Proof of Theorem \ref{thm:C-consistency}} ($C$-consistency of $\tilde{\lambda}_t$ and $\tilde{\mu}_t$).
By Theorem \ref{thm:CLT}, conditional on $\{Z_t>0\}$, the variances of $\tilde{\lambda}_t-\lambda$ and $\tilde{\mu}_t-\mu$ vanish as $t\to\infty$. This implies \eqref{C-con-lam} and \eqref{C-con-mu}, which concludes the proof.
 
\hfill$\square$ 

\bigskip

\noindent \textbf{Proof of Corollary \ref{cor-asym-var}.}
From the proof of Theorem \ref{thm:CLT}, we see that the asymptotic variances of $\tilde\lambda_t$ and $\tilde\mu_t$ (conditional on $Z_t > 0$) are, respectively, 
\[ 
  \frac{\lambda}{t(\pi^\uparrow + m - 1)} \quad \text{ and }\quad \frac{\mu}{t(\pi^\uparrow - 1)}.
\]
Because $m \geq 2$ and $\rho<0$, we have $\lambda \leq \lambda(m - 1) < \mu$. In addition, $\pi^\uparrow + m - 1 > \pi^\uparrow - 1$. We can then conclude that the asymptotic variance of $\tilde\lambda_t$ is strictly smaller than that of $\tilde\mu_t$.
\hfill $\square$

\subsection{Proofs of the results in Section \ref{sec:pk}}

\noindent\textbf{Proof of Theorem \ref{thm:C-consistency_pk}} ($C$-consistency and asymptotic normality of $\tilde{p}_{k,t}$).
For $C$-consistency, we only need to show that \begin{equation}\label{toshow}\dfrac{b^{(t)}_{k,t}}{t}\xrightarrow{p} \lambda \,p_k\,(\pi^\uparrow+k-1), \quad t \to \infty.
\end{equation}
The result then follows from 
Corollary \ref{lem:tautt}, and the Continuous Mapping Theorem.

To show \eqref{toshow}, we show that, as $t \to \infty$,
\begin{itemize}\item[\emph{(a)}] $b^{\uparrow}_{k,t}/t\xrightarrow{p} \lambda \,p_k\,(\pi^\uparrow+k-1) $, and
\item[\emph{(b)}] $\left(b_{k,t}^{(t)} - b_{k,t}^\uparrow\right)/t^p\xrightarrow{p} 0$ for any $p>0$.
\end{itemize}
To show \emph{(a)}, we use the spine decomposition\begin{equation*}
\dfrac{b_{k,t}^\uparrow}{t} = \frac{1}{t}\sum_{i=1}^{N_t} \left(\mathds{1}_{\{\tilde{\xi}_i=k\}}+\sum_{j=1}^{\tilde{\xi}_i - 1} b_{k}^{(i, j)}\right) - \frac{1}{t}\sum_{\ell=1}^{Z_t^\uparrow-1} \tilde b_{t, k}^{(\ell)},
\end{equation*}
where \( N_t \) is the number of birth events along the spine up to time \( t \),
 \( \tilde{\xi}_i \) is the number of offspring 
 generated in the \( i \)-th birth event along the spine, \( b_{k}^{(i, j)} \) is the total number of birth events until extinction that generate $k$ offspring in the \( j \)-th subcritical process born in the $i$-th birth event along the spine, and $\tilde b_{t, k}^{(\ell)}$ is the total number of birth events (until extinction) that generate $k$ offspring in the process initiated from the $\ell$-th non-spine particle alive at time $t$. The proof then follows in a similar manner to that of Lemma \ref{lem:uparrow_conv}\emph{(ii)}, noting 
  that $\mbE[\mathds{1}_{\{\tilde{\xi}_i=k\}}]=k\, p_k/m$, and $\mbE[\tilde b_{t, k}^{(\ell)}]=-\lambda\,p_k/\rho$.
  
  \smallskip
  
  To show \textit{(b)}, we observe that, for all $k\geq 2$, $b_{k,t}^{(t)}\leq b_{t}^{(t)}$ and $b_{k,t}^{\uparrow}\leq b_{t}^{\uparrow}$, Therefore, on $A_{1, t},$ 
$A_{2, q, t},$ and $
A_{3, q, t}$ with $0<q<\min\{p, |\rho|C\}$, the bounds obtained in the proof of Lemma \ref{lem:diff_t_uparrow}\emph{(ii)} apply to give
  \begin{equation} \mbE_i\left[\left|b_{k,t}^{(t)}-b_{k,t}^{\uparrow}\right|\right]\leq \mbE_i\left[B^{(t)}_{[t-C(i,q)\log t,t]}\right]+\mbE_i\left[B^\uparrow_{[t-C(i,q)\log t,t]}\right] = O(t^q\log t), \quad t \to \infty.
  \end{equation}
  \medskip
  
  For the asymptotic normality part, we again follow similar steps to the analogous result for $\tilde\lambda_t$. In this case, Corollary \ref{lem:tautt} and Slutsky's Theorem imply that \eqref{CLT-pk} will follow if we can show that
  \begin{equation}\label{42}
    \hat{Y}_{k,t}^{(t)}:=  \frac{b_{k,t}^{(t)} -p_k\, \lambda\, (\tau_t^{(t)}+(k-1)\,t)}{\sqrt{t}\sqrt{\lambda(\pi^\uparrow + k - 1)}} \xrightarrow{d} Y_k, \quad t \to \infty,
  \end{equation}
  where $Y_k\sim \mathcal{N}(0,p_k)$.
  We define the analogue of $\hat{Y}_{k,t}^{(t)}$ in the $Q$-process:
 $$ \hat{Y}_{k,t}^{\uparrow}:=  \frac{b_{k,t}^{\uparrow} -p_k\, \lambda\, [\tau_t^{\uparrow}+(k-1)\,t]}{\sqrt{t}\sqrt{\lambda(\pi^\uparrow + k - 1)}},$$
 and proceed by showing that, as $t \to \infty$, 
\begin{itemize}\item[ \emph{(c)}] $\hat{Y}^{(t)}_{k,t}-\hat{Y}_{k,t}^{\uparrow}\xrightarrow{p} 0$, and \item[\emph{(d)}] $\hat{Y}_{k,t}^{\uparrow}\xrightarrow{d} Y_k.$ \end{itemize}
Equation \eqref{42} then follows from Slutsky's Theorem.

\smallskip

To prove \emph{(c)}, we write
$$\hat{Y}^{(t)}_{k,t}-\hat{Y}^{\uparrow}_{k,t}=\dfrac{(b_{k,t}^{(t)}-b_{k,t}^{\uparrow})+p_k\,\lambda\,(\tau_t^{\uparrow}-\tau_t^{(t)})}{\sqrt{t}\sqrt{\lambda(\pi^\uparrow + k - 1)}},$$and 
we apply \emph{(b)} and Lemma \ref{lem:diff_t_uparrow}\emph{(i)} with $p=1/2$.

\smallskip

To prove \emph{(d)}, we note that the $Q$-process $Z^{\uparrow}$ is a birth-and-death process whose rate of birth associated with a jump of size $\ell-1$, $\ell\geq 2$ at population size $r$, and death rate at population size $r$, take the respective forms 
\begin{equation}\label{haha}
\lambda_r^{(\ell)} = \lambda\,p_\ell\,(r+\ell-1), \quad\text{and}\quad \mu_r = \mu\,(r-1).
\end{equation}
Assuming that a trajectory of $Z^{\uparrow}$ is continuously observed during the interval $[0,t]$, and the observations are recorded in $\vc X$, the log-likelihood is given by
\begin{align}
\ell(\vc X,p_k)&:=-\sum_{r\geq 1}\sum_{\ell\geq 2} \lambda\,p_\ell\,(r+\ell-1)\,\nu_{r, t}^\uparrow -\sum_{r\geq 2} \mu\,(r-1)\,\nu_{r, t}^\uparrow \notag \\
&\quad + \sum_{r\geq 1}\sum_{\ell\geq 2} \beta_{r,t,\ell}^\uparrow\,\log[\lambda\,p_\ell\,(r+\ell-1)]+\sum_{r\geq 2} \delta_{r,t}^\uparrow \,\log[\mu\,(r-1)],\label{LL}
\end{align}
where $\nu_{r,t}^\uparrow$ denotes the total time (during $[0,t]$) spent in state $r$, $\beta_{r,t\ell}^\uparrow$ denotes the total number of births during $[0, t]$ associated with a jump of size $\ell-1$ while in state $r$, and $\delta_{r, t}^\uparrow$ denotes the total number of deaths while in state $r$ during $[0, t]$. The MLE $\tilde{p}_{k,t}^{\uparrow}$ for $p_k$ in $Z^{\uparrow}$ is the solution of $0=\partial \ell(\vc X,p_k)/\partial p_k$, that is, 
$$0=-\sum_{r\geq 1} \lambda\,(r+k-1)\,\nu_{r,t}^\uparrow +\sum_{r\geq 1}  \dfrac{\beta_{r,t,k}^\uparrow}{p_k}.$$
Using the fact that $\sum_{r\geq 1} r\,\nu_{r, t}^\uparrow=\tau_t^\uparrow$, $\sum_{r\geq 1} \nu_{r,t}^\uparrow= t$, and $\sum_{r\geq 1}  \beta_{r,t,k}^\uparrow=b_{k,t}^\uparrow$, we obtain
\begin{equation}\label{hihi}
\tilde{p}_{k,t}^{\uparrow}=\dfrac{b_{k,t}^\uparrow}{\lambda\,[\tau_t^\uparrow+(k-1)\,t]},\end{equation}which is the analogue of $\tilde{p}_{k,t}$ in the $Q$-process.
The (expected) information function is given by $\mathcal{I}(p_k)=-\mbE[\partial^2 \ell(\vc X,p_k)/(\partial p_k)^2]=\mbE[b_{k,t}^\uparrow]/p_k^2$.
The asymptotic variance of $\tilde{p}_{k,t}^{\uparrow}$, for large $t$, is then given by $$\textrm{Var}\left(\tilde{p}_{k,t}^{\uparrow}\right)\sim \dfrac{p_k^2}{\mbE[b_{k,t}^{\uparrow}]}\sim \dfrac{p_k^2}{t\,\sum_{r\geq 1 } \lambda \, p_k\, (r+k-1)\,\mbP[Z^{\uparrow}_\infty=r]}=\dfrac{p_k}{t\,\lambda (\pi^\uparrow+k-1)}.$$ 
By Billingsley \cite[Theorem 7.3]{billingsley1961statistical},
we then have 
$$\sqrt{t\lambda(\pi^\uparrow + k-1)}\left(\tilde{p}_{k,t}^{\uparrow}-p_k\right) \xrightarrow{d} \mathcal{N}(0,p_k), \quad t \to \infty,$$
or, equivalently,
$$\sqrt{t\,\lambda (\pi^\uparrow+k-1)}\,\left\{	\dfrac{b_{k,t}^{\uparrow}-p_k\,\lambda\,[\tau_t^\uparrow+(k-1)\,t]}{\lambda\,[\tau_t^\uparrow+(k-1)\,t]}\right\}\xrightarrow{d} \mathcal{N}\left(0,p_k\right), \quad t \to \infty.
$$
Multiplying and dividing the left-hand side above by $\sqrt{t\lambda(\pi^\uparrow + k-1)}$, and using Lemma \ref{lem:uparrow_conv}\emph{(i)} together with Slutsky's Theorem, implies \emph{(d)}.

  \hfill $\square$
  
  \bigskip
  
\noindent\textbf{Proof of Corollary \ref{prop:Qcon_pk}} ($Q$-consistency of \( \hat{p}_{k,t} \)). The result follows from the fact that $$\hat{p}_{k,t} =\tilde{p}_{k,t}\, \lambda\,\dfrac{\tau_t+(k-1)\,t}{b_t}=\tilde{p}_{k,t}\, \lambda\,\dfrac{\frac{\tau_t}{t}+(k-1)}{\frac{b_t}{t}},$$ that conditional on $Z_t>0$, $\tilde{p}_{k,t}\xrightarrow{p} p_k$ (Theorem \ref{thm:C-consistency_pk}), 
$\tau_t/t \xrightarrow{p}\pi^\uparrow$ 
(Corollary \ref{lem:tautt}), and $b_t/t \xrightarrow{p}\lambda(\pi^\uparrow+m-1)$ (Lemma 
\ref{lem:uparrow_conv}\textit{(ii)}, and Lemma \ref{lem:diff_t_uparrow}\textit{(ii)} with $p=1$), and from the Continuous Mapping Theorem. \hfill $\square$

\appendix
\section{Appendix}\label{sec:appendix}
\subsection{Perron Frobenius decomposition}\label{PFdec}
Here we prove the existence of the Perron--Frobenius triple $(\rho_*, \vc u, \vc v)$ in \eqref{eigentriple} and \eqref{PF-asymp}, using an approach similar to \cite[Example 8]{CVGenCrit}.

\medskip

Let $\mathcal{L}$ denote the generator of $Z$, defined for functions $f:\NN \to \mathbb{R}$ satisfying $f(0) = 0$, by
\[
\mathcal{L}f(x) = \mu x \big(f(x-1)-f(x)\big) + \lambda x \sum_{k \ge 2} p_k \big(f(x-1+k)-f(x)\big), \qquad x \in \NN.
\]
Defining $V: \mathbb N \to \mathbb R : x \mapsto (x + 1)^\alpha$, $\alpha > 2$,  we have
\[
\mathcal L V(x) = \mu x[x^\alpha - (x + 1)^\alpha] + \lambda x \sum_{k \ge 2}p_k[(x+k)^\alpha - (x + 1)^\alpha].
\]
Now, for $x$ sufficiently large, Taylor's Theorem gives 
\[
  (x + r)^\alpha = x^\alpha + \alpha r x^{\alpha - 1} + O(r^2x^{\alpha - 2}).
\]
Hence, as $x \to \infty$, we have
\begin{align*}
 \mathcal L V(x) 
 &= \mu x[-\alpha x^{\alpha - 1} + O(x^{\alpha - 2})] +  \lambda x \sum_{k \ge 2}p_k[\alpha(k-1)x^{\alpha - 1} + O(k^2x^{\alpha - 2})] \\
 &= -\mu \alpha x^\alpha + \lambda(m-1)\alpha x^{\alpha} + O(x^{\alpha - 1}),
\end{align*}
where we have used the fact that the offspring distribution has finite variance. Since $\lambda(m-1) - \mu < 0$, it follows that there exists $M, c, c' > 0$ such that for all $x > M$
\[
  \mathcal L V(x) \le - cx^\alpha \le -c'V(x).
\]
Taking $D_0 = \{1, \dots, M\}$, it follows that there exists a constant $K > 0$ such that
\[
  \mathcal L V(x) \le -c'V(x) + K\mathbf 1_{D_0}(x).
\]

\smallskip

Theorem~5.1 of \cite{CVGenCrit} therefore applies and guarantees the existence of positive vectors $\vc u, \vc v$ and a constant $\rho_* \in \mathbb{R}$ such that 
\begin{equation}
  \vc u^\top Q = \rho_* \vc u^\top, 
  \qquad Q \vc v = \rho_* \vc v, 
  \qquad \vc u^\top \vc 1 = 1, 
  \qquad \vc u^\top \vc v = 1,
\end{equation}
where $\vc 1$ denotes the vector of ones, and $Q$ is the generator matrix of $Z$. Moreover, for each $t \ge 0$, there exists a matrix $R_t$ such that
\begin{equation}\label{PF}
  P(t) = \ee^{\rho_* t}\, \vc v \vc u^\top + R_t,
\end{equation}
and 
\begin{equation}
\|\vc e_i^\top R_t V\| \le K\, \ee^{(\rho_* - \eps)t} V(x), 
\qquad \text{for all } x \in \NN,
\end{equation}
for some $\eps > 0$, each basis vector $\vc e_i$, and a finite constant $K<\infty$.  

Finally, defining the weighted norm $\|f\|_{p,V} := \|f/V\|_p$, $p \in [1,\infty]$, we obtain
\begin{equation}\label{eq:R-bound}
\|R_t f\|_{\infty, V} \le \|f\|_{\infty, V}\, K \ee^{(\rho_*-\eps)t} 
\;\;\Leftrightarrow\;\;
\ee^{-\rho_* t}\|R_t f\|_{\infty, V} \le \|f\|_{\infty, V}\, K \ee^{-\eps t},
\end{equation}
with a similar bound for the dual norm $\|\cdot\|_{1,V}$ involving $R_t^\top$.

\subsection{MEXIT couplings between the $Q$-process and the conditioned process}\label{sec:coupling}

As explained in Section~\ref{sec:cond}, for $t \ge 0$ the conditioned process 
\[
  (Z_s^{(t)})_{0 \le s \le t} := (Z_s \mid Z_t > 0)_{0 \le s \le t}
\]
is generally difficult to analyse directly, since it evolves as a time-\emph{inhomogeneous} Markov process: its transition rates depend explicitly on the remaining time to $t$. By contrast, the $Q$-process $(Z_s^\uparrow)_{s \ge 0}$ is time-\emph{homogeneous} and therefore much more tractable. To approximate the behaviour of $(Z_s^{(t)})$ for large~$t$, we construct couplings between the two processes such that their sample paths coincide for as long as possible before eventually diverging. These are \emph{MEXIT couplings} (maximal exit time couplings), introduced in \cite{Mexit}, which we describe below.

\smallskip

Let $J_k^\uparrow$ denote the time of the $k$th jump of the $Q$-process and define the number of jumps up to time $s \ge 0$ by
\[
  K_s^\uparrow := \sup\{k \ge 0 : J_k^\uparrow \le s\}.
\]
Since we consider linear birth-and-death processes, $K_s^\uparrow$ is almost surely finite for every finite~$s$. For an initial state $i \in \NN$, we denote by 
\[
  p_i^{(s,\uparrow)}(\dd \vc{u}, \vc{x}, k) 
  := \mbP_i\!\left(K_s^\uparrow = k,\, J_\ell^\uparrow \in \dd u_\ell,\, Z_{u_\ell}^\uparrow \in x_\ell,\, 1 \le \ell \le k\right)
\]
the joint distribution of the jump times and jump sizes of the $Q$-process up to time $s$. We define analogous quantities $J_k^{(t)}$, $K_s^{(t)}$, and $p_i^{(s,t)}(\dd \vc{u}, \vc{x}, k)$ for the conditioned process $(Z_s^{(t)})$.

\smallskip

A \emph{coupling} of $(Z_s^{(t)})_{0\leq s \le t}$ and $(Z_s^\uparrow)_{0\leq s \le t}$ is a joint process 
\[
  \big(\widehat{Z}_s^{(t)}, \widehat{Z}_s^\uparrow\big)_{0 \le s \le t}
\]
with law $\widehat{\mbP}_i^{(t,\uparrow)}$ such that both marginals have the correct distributions: for every $s \le t$, the trajectory of $(\widehat{Z}_u^{(t)})_{u \le s}$ has law $p_i^{(s,t)}$, and the trajectory of $(\widehat{Z}_u^\uparrow)_{u \le s}$ has law $p_i^{(s,\uparrow)}$.  

We define the \emph{uncoupling time} as
\begin{equation}\label{uncoupling-time}
  \zeta_t := \inf\{s \le t : \widehat{Z}_s^{(t)} \neq \widehat{Z}_s^\uparrow\},
\end{equation}
with the convention that $\zeta_t = \infty$ if the two processes remain identical up to time~$t$. Thus $\zeta_t$ measures the first time at which the sample paths diverge.

The idea of a MEXIT coupling is to construct a joint law of $\big(\widehat{Z}_s^{(t)}, \widehat{Z}_s^\uparrow\big)$ that maximises the probability of agreement, or equivalently, that maximises the uncoupling time $\zeta_t$. In this sense, under the MEXIT coupling the random variable $\zeta_t$ stochastically dominates the uncoupling time under any other coupling. 

\begin{figure}[t]
\centering\includegraphics[width=0.7\textwidth]{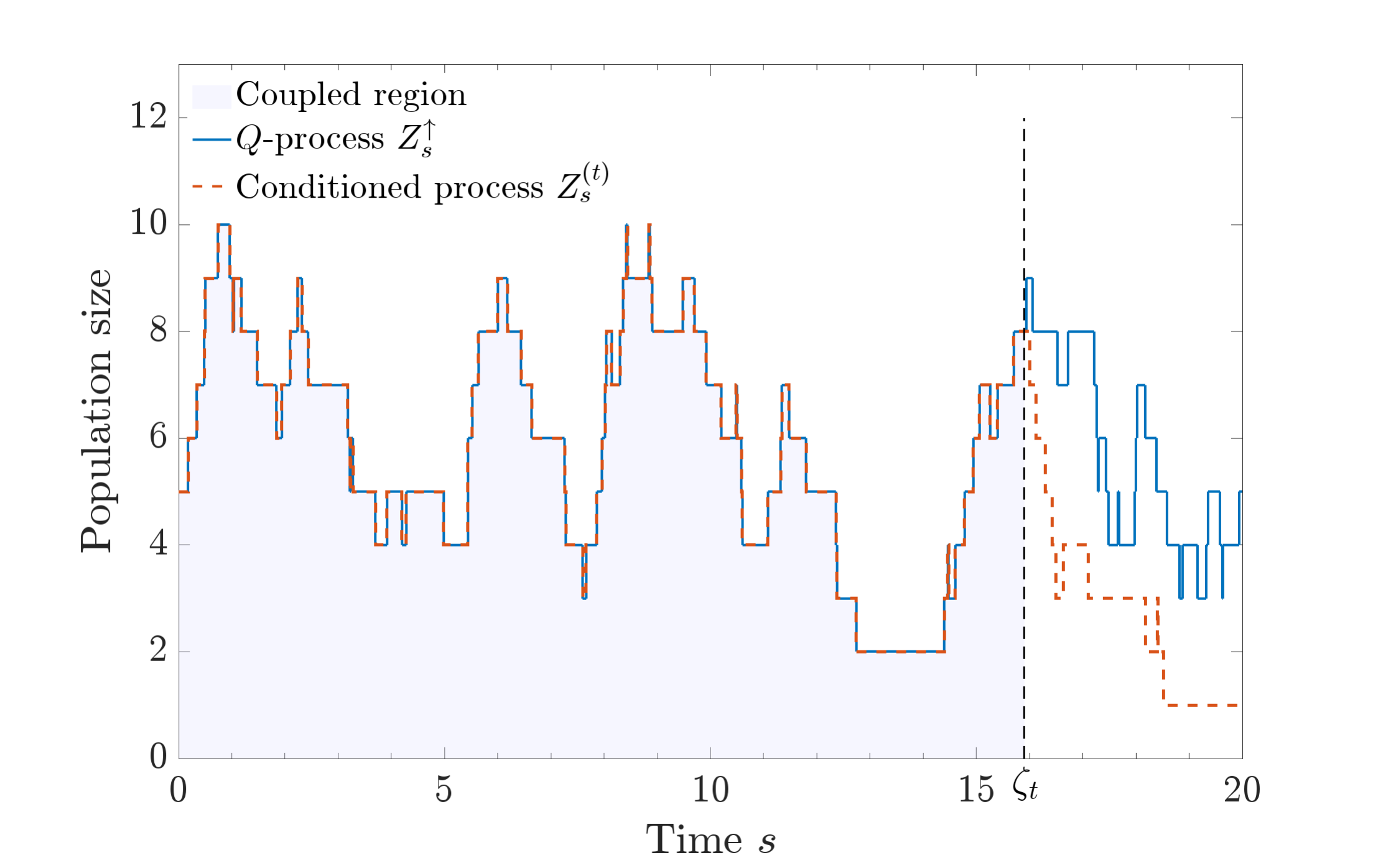}
\caption{\label{f2}\textbf{MEXIT coupling:} $Q$-process $Z_s^\uparrow$ and conditioned process $Z_s^{(t)}$ for $t=20$. The uncoupling time $\zeta_t$ is the first divergence between the two trajectories.}
\end{figure}

By \cite{Mexit}, for each $t\geq 0$, there exists such a MEXIT coupling of $(Z_s^{(t)})_{s \le t}$ and $(Z_s^\uparrow)_{s \le t}$, which we denote by $\widehat{\mbP}_i^{(t,\uparrow)}$. Since MEXIT couplings are maximal, the probability of uncoupling by time $s$ coincides with the total variation distance between the trajectory distributions of the two processes up to $s$, namely
\begin{equation}
  \widehat{\mbP}_i^{(t, \uparrow)}(\zeta_t \le s) 
  = \frac12 \sum_{k = 0}^\infty
   \int_0^s \dd u_1 \int_{u_1}^s \dd u_2 \cdots \int_{u_{k-1}}^s \dd u_k 
   \sum_{\vc{x} \in \NN^k}\!\big|p_i^{(s, t)}(\dd \vc{u}, \vc{x}, k) - p_i^{(s, \uparrow)}(\dd \vc{u}, \vc{x}, k)\big|.
  \label{mexit}
\end{equation}
This characterisation shows that the MEXIT coupling keeps the two processes identical for as long as possible and provides a way to quantify how close $(Z_s^{(t)})$ is to $(Z_s^\uparrow)$ when $t$ is large. Figure~\ref{f2} illustrates a MEXIT coupling.

\smallskip

The next result provides quantitative control on the uncoupling time.

\begin{prop}\label{thm:mexit}
For each $i \in \NN$ the MEXIT couplings $\widehat{\mbP}_i^{(t,\uparrow)}$ of $(Z_s^{(t)})_{s \le t}$ and $(Z_s^\uparrow)_{s \le t}$ , $t \ge 0$, are such that:

\begin{enumerate}
  \item[(i)] For every $\ell \ge 0$, the probability that uncoupling occurs before time $t-\ell$ satisfies
  \begin{equation}
    \lim_{t \to \infty}\widehat{\mbP}_i^{(t,\uparrow)}(\zeta_t \le t-\ell) 
    = \frac{\ee^{-\rho \ell}}{2} \sum_{j = 1}^{\infty} u_j \big|\vc{e}_j^T R_\ell \vc{1}\big|.
    \label{uncoupling-bound}
  \end{equation}
  
  \item[(ii)] The probability of ever uncoupling converges to
  \begin{equation}
    \lim_{t \to \infty}\widehat{\mbP}_i^{(t,\uparrow)}(\zeta_t < \infty) 
    = \frac12 \sum_{j = 1}^\infty u_j \,|1 - v_j|.
    \label{uncoupling-finite}
  \end{equation}

  \item[(iii)] For each $q > 0$, there exist constants $C(i, q)$ and $T(i, q)$ such that
  \begin{equation}
    \widehat{\mbP}_i^{(t,\uparrow)}\!\left(\zeta_t \le t - C(i, q)\log t \right) \le \frac{1}{t^q},
    \qquad t \ge T(i, q).
    \label{uncoupling-good}
  \end{equation}
\end{enumerate}
\end{prop}

\smallskip

In words, this proposition shows that one can construct couplings in which the conditioned process and the $Q$-process remain identical with high probability over most of the time interval $[0,t]$. In particular, part~\textit{(iii)} implies that, with probability tending to one, the processes only uncouple in the final $C(i, q) \log t$ units of time before~$t$, so that for large~$t$ the $Q$-process provides an accurate and tractable approximation of the conditioned process. Its proof follows similar ideas as used in~\cite{BHM22}.

\medskip

\noindent\textbf{Proof.}

\noindent\emph{Step 1. Path decomposition.}  
We consider the right-hand side of \eqref{mexit}, and decompose $p_i^{(s, t)}$ into
\[
  p_i^{(s, t)}(\dd \vc{u}, \vc{x}, k) = 
  \frac{P_{i, x_1}(u_1)P_{x_1, x_2}(u_2-u_1)\cdots P_{x_{k -1}, x_k}(u_k-u_{k-1})\ee^{-x_k(\lambda + \mu)(s-u_k)}
   \vc{e}_{x_k}^TP(t-s)\vc{1}}{\vc{e}_i^\top P(t)\vc{1}},
\]
Using a similar decomposition for $p_i^{(s, \uparrow)}$, we obtain the same expression as for $p_i^{(s, t)}$ with $\vc{e}_{x_k}^TP(t-s)\vc{1}/(\vc{e}_i^\top P(t)\vc{1})$ replaced by $v_{x_k}/(\ee^{\rho s}v_i)$.

Substituting this into~\eqref{mexit} and simplifying yields
\begin{align}
  \widehat{\mbP}_i^{(t,\uparrow)}(\zeta_t \le s)
  &= \frac12 \sum_{j\ge1} \vc{e}_i^\top P(s)\vc{e}_j
  \left|\frac{\vc{e}_j^\top P(t-s)\vc{1}}{\vc{e}_i^\top P(t)\vc{1}}
  - \frac{v_j}{\ee^{\rho s} v_i}\right|\notag\\
  &= \frac12 \sum_{j\ge1} \Big(\ee^{\rho s}u_j v_i + \vc{e}_i^\top R_s \vc{e}_j\Big)
  \left|\frac{\ee^{\rho s} v_i \vc{e}_j^\top R_{t-s}\vc{1}
  - v_j \vc{e}_i^\top R_t\vc{1}}
  {\ee^{\rho s} v_i\big(v_i \ee^{\rho t} + \vc{e}_i^\top R_t\vc{1}\big)}\right|,
  \label{eq:step1}
\end{align}
where we used the decomposition~\eqref{PF}, recalling that $\rho_*=\rho$.

\smallskip

\noindent\emph{Step 2. Upper and lower bounds on $\widehat{\mbP}_i^{(t,\uparrow)}(\zeta_t \le s)$.}  
Our aim is to find upper and lower bounds for the right-hand side of \eqref{eq:step1} that both converge as $t\to\infty$ to the right-hand side of~\eqref{uncoupling-bound}. 

Starting with the upper bound, using H\"older's inequality and the bound on $R_s$ from~\eqref{eq:step1}, we have
\begin{equation}
 |\vc{e}_i^\top R_s \vc{e}_j| 
 \le \|\vc{e}_j^\top\|_{\infty, V} \|R_s^\top \vc{e}_i\|_{1, V}
 \le K\, \frac{V(i)}{V(j)}\,\ee^{(\rho - \eps) s},
\label{R-bound-1}
\qquad |\vc{e}_i^\top R_s \vc{1}| \le K\, V(i) \,\ee^{(\rho-\eps) s}
\end{equation}
for some $K<\infty$ and $\eps>0$.
Now, note that for any $\eta > 0$, we can choose $T(i)$ sufficiently large such that for all $t > T(i)$, 
\[
  V(i)K\ee^{-\eps t} < \eta v_i.
\]
Now fix $\ell = t - s$. For $t > T(i) + \ell$ (which implies that $s > T(i)$), the right-hand side of~\eqref{eq:step1} is bounded above by
\begin{align}
\frac12 \sum_{j = 1}^\infty&\left(\ee^{\rho s}u_jv_i + \frac{\eta\, v_i\,\ee^{\rho s}}{V(j)}\right)\left(\frac{\ee^{\rho s}v_i |\vc{e}_j^\top R_\ell \mathbf{1}| + v_jv_i \ee^{\rho t}\eta}{\ee^{\rho s}v_i(v_i\ee^{\rho t} - v_i \ee^{\rho t}\eta)} \right)\notag\\
&= \frac12 \sum_{j = 1}^\infty\left(u_j + \frac{\eta}{V(j)}\right)\left(\frac{\ee^{-\rho \ell}|\vc{e}_j^\top R_\ell \mathbf{1}| + v_j \eta}{1 - \eta} \right)\notag\\
&= \frac{\ee^{-\rho \ell}}{2(1 - \eta)}\sum_{j = 1}^\infty u_j |\vc{e}_j^\top R_\ell \mathbf{1}| 
+ \frac{\eta}{2(1 - \eta)} 
+ \frac{\eta\ee^{-\rho \ell}}{2(1 - \eta)}\sum_{j = 1}^\infty \frac{|\vc{e}_j^\top R_\ell \mathbf{1}|}{V(j)} 
+ \frac{\eta^2}{2(1 - \eta)}\sum_{j = 1}^\infty \frac{v_j}{V(j)},
\label{UB}
\end{align}
where we have used the fact that $\vc u^T\vc v = 1$ to obtain the second term in the final step.

A similar argument with the reverse triangle inequality gives the corresponding lower bound to \eqref{eq:step1}:
\begin{align}
\frac12 \sum_{j = 1}^\infty &\left(u_j - \frac{\eta}{V(j)}\right)\left(\frac{\ee^{-\rho \ell}|\vc{e}_j^\top R_\ell \mathbf{1}| - v_j\eta}{1 + \eta} \right)\notag \\
 &= \frac{\ee^{-\rho \ell}}{2(1 - \eta)}\sum_{j = 1}^\infty u_j |\vc{e}_j^\top R_\ell \mathbf{1}| 
 - \frac{\eta}{2(1 - \eta)} 
 - \frac{\eta\ee^{-\rho \ell}}{2(1 - \eta)}\sum_{j = 1}^\infty \frac{|\vc{e}_j^\top R_\ell \mathbf{1}|}{V(j)} 
 + \frac{\eta^2}{2(1 - \eta)}\sum_{j = 1}^\infty \frac{v_j}{V(j)}.
\label{LB}
\end{align}

\smallskip

\noindent\emph{Step 3. Finiteness of error terms.}
If we can show the third and fourth terms on the right-hand side of~\eqref{UB} and~\eqref{LB} are finite, the result follows by taking $\eta$ arbitrarily small. 

Let us first consider the third term. From~\eqref{PF}, we have that for any $\ell \ge 1$
\begin{equation}
  R_\ell = P(\ell) - \ee^{\rho \ell}\vc v \vc u^T, \qquad 
  |\vc{e}_j^T P(\ell) \mathbf{1}| \le 1, \qquad \text{and} \qquad 
  \vc{e}_j(\vc v\vc u^T)\mathbf{1} = v_j, 
\label{3ids}
\end{equation}
so
\begin{equation}
\sum_{j = 1}^\infty \left|\frac{\vc{e}_j^\top R_\ell \mathbf{1}}{V(j)}\right| 
= \sum_{j = 1}^\infty \left|\frac{\vc{e}_j^\top (P(\ell) - \ee^{\rho \ell}\vc v\vc u^\top)\mathbf{1}}{V(j)}\right|
\le \sum_{j = 1}^\infty \frac{1 + \ee^{\rho \ell}v_j}{V(j)}<\infty
\label{third-term}
\end{equation}
since $v_j = c\, j$ for some finite constant $c$ and $V(j) = (1 + j)^\alpha$ with $\alpha > 2$. Finiteness of the fourth term follows similarly. 

\smallskip

\noindent\emph{Step 4. Limits.}  
Using \eqref{mexit} with $s=t-\ell$ and letting $\eta\to0$ in~\eqref{UB} and~\eqref{LB} yields
\[
  \lim_{t\to\infty}\widehat{\mbP}_i^{(t,\uparrow)}(\zeta_t \le t-\ell)
  = \frac{\ee^{-\rho \ell}}{2}\sum_{j\ge1} u_j |\vc{e}_j^\top R_\ell \vc{1}|,
\]
which is claim (i). We obtain (ii) by choosing $\ell=0$.

Finally, for part (iii), we take $t - s = \ell = C(i, q)\log t$ and note that for any $q, \gamma > 0$, we may choose $C(i, q)$ sufficiently large so that
\begin{equation}
V(i)K\ee^{-\eps \lfloor C(i, q)\log t \rfloor} \le \frac{\gamma v_i}{t^q},
\label{3bound1}
\end{equation}
for $t$ sufficiently large. This then implies that 
\begin{equation}
V(i)K\ee^{-\eps s} \le \frac{\gamma v_i}{t^q} \qquad \text{and} \qquad V(i)K\ee^{-\eps t} \le \frac{\gamma v_i}{t^q},
\label{3bound2}
\end{equation}
for $t$ sufficiently large. 

Combining~\eqref{eq:step1},~\eqref{R-bound-1},~\eqref{3bound1} and~\eqref{3bound2}, we have
\begin{align}
  \widehat{\mbP}_i^{(t, \uparrow)}(\zeta_t \le s) 
  &\le \frac12 \sum_{j = 1}^\infty \left(\ee^{\rho s}u_j v_i + \frac{\gamma t^{-q}v_i \ee^{\rho s}}{V(j)}\right)\left(\frac{\ee^{\rho t} v_i v_j \gamma t^{-q} + v_jv_i\ee^{\rho t}\gamma t^{-q}}{\ee^{\rho s}v_i (v_i \ee^{\rho t} - v_i \ee^{\rho t}\gamma t^{-q})} \right)\notag \\
  &\le \frac{\gamma t^{-q}}{1 - \gamma t^{-q}}\left(\sum_{j = 1}^\infty u_j v_j + \gamma t^{-q}\sum_{j = 1}^\infty \frac{v_j}{V(j)}\right).\label{3bound3}
\end{align}
Again, using the fact that $\vc u^T\vc v = 1$ along with similar arguments used to obtain finiteness of~\eqref{third-term},~\eqref{3bound3} is bounded above by $K_2\gamma t^{-q}$ for some constant $K_2$ and $t$ sufficiently large. Choosing $\gamma \le K_2^{-1}$ yields the result. \hfill $\square$

\section*{Acknowledgements}
Both authors would like to thank the Australian Research Council (ARC) for support through the Discovery Project DP200101281 and PHC FASIC chercheurs 2022 programme, project number 48421UH.

\bibliography{ref}{}
\bibliographystyle{plain}

\end{document}